\documentclass{article}

\usepackage{arxiv}

\usepackage[utf8]{inputenc} 
\usepackage[T1]{fontenc}    
\usepackage{hyperref}       
\usepackage{url}            
\usepackage{booktabs}       
\usepackage{amsfonts}       
\usepackage{nicefrac}       
\usepackage{microtype}      
\usepackage{lipsum}		
\usepackage{amsthm}
\usepackage{graphicx}
\usepackage{doi}
\usepackage{physics}
\usepackage{amssymb}
\usepackage{pifont}

\newcommand{\pd}[2]{\frac{\partial #1}{\partial #2}}
\newcommand{\pdd}[2]{\frac{\partial^2 #1}{\partial #2^2}}

\usepackage{setspace}
\title{Extending the explicit constraint force method to inverse problems}

\date{} 					

\author{
    Conor Rowan \\
	Aerospace Engineering\\
	University of Colorado Boulder\\
    3775 Discovery Drive \\
	Boulder, CO 80309 \\
	\texttt{conor.rowan@colorado.edu} \\
}

\renewcommand{\headeright}{}
\renewcommand{\undertitle}{}
\renewcommand{\shorttitle}{}

\begin{document}
\maketitle

\begin{abstract}
    Recently, the explicit constraint force method (ECFM) was introduced as a principled approach to solution reconstruction in the presence of missing physics. In solution reconstruction, parameters of a physical model are estimated from sparse measurement data as a means to obtain the full solution field. In contrast, inverse problems target the missing parameters and estimate the solution along the way. Noting the similarity of the mathematical formulations of these two tasks, we investigate the use of ECFM to solve inverse problems. First, we compare the ECFM formulation of the inverse problem to a standard approach using two numerical examples. The first example provides an extension of ECFM to dynamic problems, and the second offers a novel approach to treat noisy measurement data. Next, we introduce a method to solve inverse problems for which the parameterized model has stochastic components. This approach is based on constraint forces and the polynomial chaos expansion, and is illustrated with another numerical example. Finally, we discuss extensions of ECFM to recover missing boundary conditions and domain geometries from measurement data, which are shown to be special cases of problems treated previously in the literature. The purpose of this work is to extend the mathematical framework of ECFM to novel applications and to gauge the method's viability as an alternative strategy for inverse analysis.
\end{abstract}


\section{Introduction}

\paragraph{} Philosopher and physicist Mario Bunge defines the task of an inverse problem as identifying the ``hidden causes of observable effects'' and gives some examples of inverse problems familiar from daily life: planning, diagnosis, invention, and design \cite{bunge_inverse_2019}. According to Bunge, these kinds of problems are ubiquitous and generally harder than the corresponding forward problem. Within scientific computing, inverse problems are frequently used to estimate parameters describing phenomena in physical systems governed by ordinary or partial differential equations (PDEs). For example, an inverse problem can be solved to estimate the material properties of a linearly elastic solid using a discrete set of strain measurements \cite{hematiyan_new_2017}. Another common inverse problem from geophysics is full-wave inversion, in which the material properties of a medium are inferred from data on wave propagation \cite{mercier_designing_2025}. Recovering a missing source term from measurement data is a related inverse problem with applications in seismology and heat transfer, but is challenging due to its ill-posedness \cite{jang_solution_2000, shidfar_solving_2010}. Because all physical models depend on empirical parameters that are calibrated from measurement data, inverse problems are a mainstay of scientific computing. This realization challenges the bright line often drawn between ``physics-based'' and ``data-driven'' models in the scientific machine learning community \cite{karniadakis_pinns_2020}. From Newton's law of gravitation to the Navier-Stokes and Maxwell equations, all scientific models are data-driven in some sense, though these models typically have fewer parameters and the ``training'' process is an inverse problem rather than the machine learning approaches of modern data-driven models \cite{li_fourier_2021, lu_deeponet_2021}. 

\paragraph{} At the highest level, mathematical formulations of inverse problems split into two broad categories: Bayesian and deterministic. A typical Bayesian inverse problem uses a parameterized model, a known distribution of measurement noise, a prior distribution on the unknown parameters, and Bayes' rule to construct the posterior distribution on the parameters of interest \cite{dashti_bayesian_2015}. The Bayesian approach to inverse problems has been studied extensively in the scientific community, for example in biomechanics \cite{franck_multimodal_2017}, fracture mechanics \cite{khodadadian_bayesian_2020}, and full-wave inversion \cite{zhang_bayesian_2021}. It is desirable because uncertainty in the parameter estimate is naturally quantified, which can be an asset for downstream uncertainty quantification tasks. The presence of a prior distribution on the parameters also allows for the seamless incorporation of prior knowledge, which may be essential in obtaining parameter estimates if measurement data is scarce. Furthermore, Bayesian methods are convenient for treating ill-posed problems, as multi-modal posterior distributions avoid overconfident point estimates which are non-unique \cite{sun_local_2022}. One drawback of the Bayesian approach is that the posterior distribution does not provide insight into whether the parameterized model is ``consistent'' with the measurement data. By consistent, we mean that there is some parameter setting of the model which recovers the data up to random measurement noise. Post-processing steps involving the posterior predictive distribution are required to do this work of ``model checking'' \cite{gelman_philosophy_2013}. Though the posterior distribution avoids overconfident point estimates of the parameters, it does not necessarily guard against overconfidence about the form of the model itself. One standard approach is to introduce an additive discrepancy such as a Gaussian process to account for the model misspecification \cite{kennedy_bayesian_2001, zou_correcting_2024}. A more radical alternative is to use a Gibbs posterior distribution, which does not require knowledge of the data generation process (in other words, the likelihood) and is thus more robust to misspecified models \cite{jiang_gibbs_2008}.


\paragraph{} While we maintain our interest in misspecified/inconsistent models, we will not pursue the Bayesian approach to inverse problems any further, as the remainder of this paper is focused on deterministic formulations. In place of using Bayes' rule, deterministic inverse problems find a point estimate of the parameters by choosing them such that the predictions of the model agree as closely as possible with the measurement data \cite{vogel_computational_2002}. Deterministic inverse problems have been particularly popular in the physics-informed machine learning literature in recent years, owing to their convenient incorporation into the standard pipeline of physics-informed neural network training. The original physics-informed neural networks (PINNs) paper recovered the missing viscosity from the Navier-Stokes equations using flow data \cite{raissi_physics-informed_2019}. This was followed up by extensions to conjugate heat transfer \cite{cai_physics-informed_2021}, contact mechanics \cite{sahin_solving_2024}, and beam bending \cite{zhou_data-guided_2024}. It is typically assumed when formulating a deterministic inverse problem that the parameterized model is consistent with the measurement data, such that the mismatch between the predictions of the calibrated model and the measurement data resembles random noise. One exception to this in the PINNs literature is \cite{zou_correcting_2024}, which introduces an additive correction to account for the discrepancy between the solution of the calibrated PDE and the measurement data. 

\paragraph{} Whereas an inverse problem fits the solution of a PDE in the process of estimating the missing model parameters, a related problem called solution reconstruction does the opposite. With solution reconstruction, the goal is to obtain a full-field estimate of the state from sparse measurement data, which is often accomplished by estimating the model parameters. In other words, inverse problems seek the missing parameters and use the solution field en route, whereas solution reconstruction uses the parameters to predict the solution. Solution reconstruction has also been a topic of interest in the PINNs community. In \cite{raissi_hidden_2020}, the velocity and pressure fields are interpolated with a neural network and estimated from image data. In \cite{yeung_physics-informed_2022}, neural networks are used to infer the full-field state of the hydraulic head of groundwater flow by estimating a spatially varying constitutive parameter. Similarly, \cite{ehlers_data_2024} estimates coefficients in the governing equation as a step in reconstructing the full-field state of the velocity field of shallow water waves. Deterministic inverse problems and solution reconstruction are equivalent when the parameterized model is consistent with the measurement data. However, when the parameterized model is inconsistent with the measurement data, the recovered parameters may cease to be meaningful and are treated only as a means to reconstruct the state field. Recently, the explicit constraint force method (ECFM) was introduced as a principled approach to performing solution reconstruction in the presence of misspecified models \cite{rowan_physics-informed_2025}. ECFM introduces additional source terms into the system which enforce that the predictions of the parameterized model agree with the measurement data. The model parameters are then chosen such that the magnitude of the constraint force is minimal, and the value of the constraint force at the minimum provides a measure of consistency between the parameterized system and the data. It was also shown that this formulation of the solution reconstruction problem allowed for the use of the variational energy objective, unlike traditional approaches. 

\paragraph{} While ECFM was first introduced for solution reconstruction, the close resemblance of solution reconstruction to deterministic inverse problems raises questions regarding the method's utility in this unexplored context. As such, the purpose of this paper is two-fold:

\begin{enumerate}
    \item We introduce ECFM as a novel formulation of deterministic inverse problems and compare its performance to traditional methods relying on the mismatch between the model predictions and the data;
    \item We discuss extensions of ECFM to dynamic problems, stochastic models, missing boundary conditions and domain geometries, as well as introduce a new approach to handling noisy measurement data.
\end{enumerate}

In \cite{rowan_physics-informed_2025}, only the recovery of parameterized differential operators and source terms in the static setting was considered, thus these extensions represent an expansion of the scope of problems that ECFM can treat. The rest of this paper is organized as follows. In Section 2, we introduce an abstract formulation of the standard PDE-constrained inverse problem, as well as sketch our approach based on constraint forces. In Section 3, we provide three numerical examples of inverse problems based on ECFM. In the first two, we compare our proposed method against a standard formulation based on the squared error between the model predictions with the data. In the third example, we illustrate the use of ECFM in the case that the parameterized model has stochastic components. Each numerical example represents an extension of ECFM to a new class of problems. Then, in Section 4, we discuss how recovering boundary conditions and domain geometry from measurement data reduces to the previously-treated cases of parametric differential operators and/or source terms. In Section 5, we summarize our findings and provide directions for future work.


\section{Two approaches to the inverse problem}

\paragraph{} In this section, we introduce a standard approach to deterministic inverse problems, as well as the proposed ECFM formulation. This is done abstractly, as we rely on the subsequent numerical examples to make the discussion concrete. Our purposes here are to illustrate the differences between the two methods and introduce notation that is used throughout the paper. To this end, we take the system of interest to obey the following initial and/or boundary value problem:

\begin{equation}\label{abstract_bvp}
\begin{aligned}
    \mathcal{N}( \mathbf{u})(\mathbf{y}) = \mathbf{0}, \quad \mathbf{y} \in \Omega,\\
    \mathcal{B}(\mathbf{u})(\mathbf{y}) = \mathbf{0}, \quad \mathbf{y} \in \partial \Omega,
\end{aligned}
\end{equation}

\noindent where $\mathbf{u}$ is the system state, $\mathbf{y}$ are coordinates that may involve space and/or time, $\Omega$ is the domain on which the state is defined, $\mathcal{N}(\cdot)$ is a differential operator, and $\mathcal{B}(\cdot)$ is a boundary operator which may involve initial and/or boundary conditions. We are deliberately vague about whether the system is static or dynamic, as subsequent numerical examples involve both sorts of problems. Eq. \eqref{abstract_bvp} is the governing equation for the system from which the measurement data used in the inverse problem is taken. Now, we measure the system state $\mathbf{u}(\mathbf{y})$ at $C$ distinct points, which is given by

\begin{equation}\label{measurement}
    \mathbf{V} = \{\mathbf{v}_i \}_{i=1}^C  = \Big \{\mathbf u (\mathbf y_i) + \boldsymbol \xi_i \Big \}_{i=1}^C = \mathcal M(\mathbf u(\mathbf y)) + \boldsymbol \Xi, 
\end{equation}

\noindent where $\{\mathbf{y}_i \}_{i=1}^C$ are the positions of measurements, $\mathcal M$ is a measurement functional, and $\boldsymbol \Xi=\{ \boldsymbol \xi_i\}_{i=1}^C$ is the additive measurement noise. The noise is taken to be zero-mean Gaussian with no correlations between measurement points or components of the state, i.e., $\boldsymbol \xi_i \overset{\text{i.i.d.}}{\sim} N(\mathbf{0}, \sigma^2 \mathbf{I})$ where $\sigma^2$ is the variance.

\paragraph{} In order to set up the inverse problem, we introduce an initial and/or boundary value problem parameterized by $\boldsymbol \epsilon$, where these are the parameters to be recovered. This reads 

\begin{equation}\label{parameterized_bvp}
\begin{aligned}
    \mathcal{G}( \mathbf{w} ; \boldsymbol \epsilon)(\mathbf{y}) = \mathbf{0}, \quad \mathbf{y} \in \Omega(\boldsymbol \epsilon),\\
    \mathcal{Q}(\mathbf{w} ; \boldsymbol \epsilon)(\mathbf{y}) = \mathbf{0}, \quad \mathbf{y} \in \partial \Omega(\boldsymbol \epsilon),
\end{aligned}
\end{equation}

\noindent where $\mathbf{w}$ is the state, $\mathcal{G}(\boldsymbol \cdot;\boldsymbol \epsilon)$ is the parameterized differential operator, $\mathcal{Q}(\boldsymbol \cdot ; \boldsymbol \epsilon)$ is the parameterized boundary operator, and $\Omega(\boldsymbol \epsilon)$ is the parameterized domain. We call $\boldsymbol \epsilon$ the ``model parameters'' going forward. Using Eqs. \eqref{measurement} and \eqref{parameterized_bvp}, the standard formulation of the deterministic inverse problem \cite{vogel_computational_2002} is given by

\begin{equation}\label{traditional}
\begin{aligned}
&\underset{\boldsymbol{\epsilon},\, \mathbf{w}(\mathbf y)}
    {\text{argmin }}\;\frac{1}{2}\, 
     \lVert \mathbf V - \mathcal{M}(\mathbf{w}(\mathbf y))\rVert^{2} \\[6pt]
\text{s.t.} \qquad 
& \mathcal{G}(\mathbf{w}; \boldsymbol{\epsilon})(\mathbf{y}) = 0,
    \quad \mathbf{y}\in \Omega(\boldsymbol{\epsilon}), \\[4pt]
& \mathcal{Q}(\mathbf{w}; \boldsymbol{\epsilon})(\mathbf{y}) = 0,
    \quad \mathbf{y}\in \partial\Omega(\boldsymbol{\epsilon}) ,
\end{aligned}
\end{equation}

\noindent where $\lVert \boldsymbol \cdot \rVert$ is the Frobenius norm in the case of a vector-valued state $\mathbf w$. This PDE-constrained optimization problem is interpreted as follows: find the model parameters $\boldsymbol \epsilon$ such that when the governing equation of Eq. \eqref{parameterized_bvp} is solved, the measured state agrees as closely as possible with the data $\mathbf V$. An advantage of this approach is its robustness to noise in the case that the parameterized model is consistent with the data. When there is measurement noise, the objective in Eq. \eqref{traditional} is not driven down to zero because the parameterized model is not capable of exactly reproducing the noisy measurement data. This is no cause for concern, however, as the discrepancy between the model and the data is zero mean in expectation at the solution of Eq. \eqref{traditional}. Thus, assuming consistency of the parameterized model and a sufficiently rich set of measurements, we expect to recover the true parameters both in the noisy and noise-free cases. We remark that there is no reason to expect accurate recovery of model parameters when the model is misspecified.

\paragraph{} Like the case of solution reconstruction, the formulation of the inverse problem based on constraint forces is more sensitive to noise \cite{rowan_physics-informed_2025}, and thus the noisy and noise-free cases require separate treatment. For the sake of exposition, we present the noise-free case here, and defer a discussion of noisy measurement data to a numerical example in Section 3. We remark that one approach to handling noisy data is introduced in \cite{rowan_physics-informed_2025}. This approach is rather complex from both a theoretical and implementation standpoint, as it involves implementing a custom solution procedure analogous to inequality-constrained optimization with Lagrange multipliers. Owing to its complexity, this is not pursued further in the current work, as we propose a more straightforward alternative in the following section. For now, the noise variance is taken to be zero $(\sigma^2=0)$. Following \cite{rowan_physics-informed_2025}, ECFM modifies the governing equation of Eq. \eqref{parameterized_bvp} to read

\begin{equation}\label{parameterized_bvp_ecfm}
\begin{aligned}
    \mathcal{G}( \mathbf{w} ; \boldsymbol \epsilon)(\mathbf{y}) + \sum_{i=1}^C \boldsymbol \lambda_i \Gamma(\mathbf{y} - \mathbf{y}_i) = \mathbf{0}, \quad \mathbf{y} \in \Omega(\boldsymbol \epsilon),\\
    \mathcal{Q}(\mathbf{w} ; \boldsymbol \epsilon)(\mathbf{y}) = \mathbf{0}, \quad \mathbf{y} \in \partial \Omega(\boldsymbol \epsilon),
\end{aligned}
\end{equation}

\noindent where $\Gamma(\mathbf{y}-\mathbf{y}_i)$ are ``constraint forces'' centered on the measurement locations and $\boldsymbol \lambda_i$ are the constraint force magnitudes. As an analogy to the method of Lagrange multipliers---which, when solving inverse or solution reconstruction problems, also introduces source terms that are scaled to enforce data constraints---the constraint force magnitudes are chosen such that the system state $\mathbf w$ agrees with the measurement data at any setting of the model parameters $\boldsymbol \epsilon$. In contrast with the ``soft constraint'' on the measurement data given in Eq. \eqref{traditional}, ECFM takes the measurement data as the ground truth, but does not assume that there is any setting of the model parameters that is capable of reproducing it. The optimal setting of the model parameters is then taken to be the one which minimizes the total constraint force:

\begin{equation}\label{ecfm}
\begin{aligned}
&\underset{\boldsymbol{\epsilon},\boldsymbol \lambda, \mathbf{w}(\mathbf y)}
    {\text{argmin }}\;\frac{1}{2}\, 
     \lVert \boldsymbol \lambda\rVert^{2} \\[6pt]
\text{s.t.} \qquad 
& \mathcal{G}(\mathbf{w}; \boldsymbol{\epsilon})(\mathbf{y})+ \sum_{i=1}^C \boldsymbol \lambda_i \Gamma(\mathbf{y} - \mathbf{y}_i) = 0,
    \quad \mathbf{y}\in \Omega(\boldsymbol{\epsilon}), \\[4pt]
& \mathcal{Q}(\mathbf{w}; \boldsymbol{\epsilon})(\mathbf{y}) = 0,
    \quad \mathbf{y}\in \partial\Omega(\boldsymbol{\epsilon}) ,\\[4pt]
& \mathcal M(\mathbf w(\mathbf y)) - \mathbf V = \mathbf 0.
\end{aligned}
\end{equation}

Though initially explored as a technique for solution reconstruction, Eq. \eqref{ecfm} is equally valid as a formulation for the inverse problem. In words, this optimization problem states: find the model parameters $\boldsymbol \epsilon$ which minimize the magnitude of the source terms required to make the solution field $\mathbf w$ coincide with the measurement data $\mathbf V$. To use the language of solid mechanics, whereas the standard inverse problem of Eq. \eqref{traditional} minimizes the discrepancy with the data in the sense of ``displacement,'' the constraint force formulation of Eq. \eqref{ecfm} finds the model parameters that minimize the discrepancy with the data in the sense of  ``force.'' In the opinion of the author, the overarching philosophy of inverse problems is tuning parameters to minimize a model's discrepancy with measurement data. When viewed in this light, there are many ways to define discrepancy. For example, small displacement discrepancies may require large constraint forces to correct. Similarly, depending on the system of interest and the position of measurements, large discrepancies in the displacement may require small constraint forces to correct. Is it clear that there is a singular correct definition of the discrepancy of the parameterized system with the measurement data? Certainly, the ECFM approach has yet to be explored in the context of inverse problems. The following section represents an initial exploration of inverse problems with constraint forces, and introduces three extensions of ECFM along the way.


\section{Numerical examples}

\subsection{Dynamics}

\paragraph{} The first numerical example is that of a classical time-dependent PDE. The measurement data in the inverse problem arises from a scalar velocity field $u(x,t)$, which is taken to obey Burgers' equation with homogeneous Dirichlet boundaries:

\begin{equation}\label{burgers}
    \pd{u}{t} + u \pd{u}{x} = \nu \pdd{u}{x} + s(x,t), \quad u(0,t) = u(1,t)= 0, \quad u(x,0) = u_0(x),
\end{equation}

\noindent where $[0,1]$ is the spatial domain, $\nu$ is the viscosity, $s(x,t)$ is a source term, and $u_0(x)$ is the initial velocity field. The velocity is observed without noise at spatial points $\{ x_i \}_{i=1}^C$ and continuously in time. The measurement data is thus $ v_i(t) = u( x_i,t)$. The goal of the inverse problem is to estimate the viscosity and the scaling of the source term from sparse measurements. The system is parameterized with

\begin{equation}\label{burgers_param}
    \pd{w}{t} + w \pd{w}{x} = 10^{-\epsilon_1} \pdd{w}{x} + \epsilon_2 s(x,t), \quad w(0,t) = w(1,t)= 0, \quad w(x,0) = u_0(x),
\end{equation}

\noindent where $w(x,t)$ is the velocity field used in recovering the parameters, and prior knowledge that the viscosity is small informs our choice to estimate its logarithm. The velocity field is discretized with

\begin{equation*}
    \hat w( x;\boldsymbol \theta (t) ) = \sum_{i=1}^N \theta_i(t) f_i(x),
\end{equation*}

\noindent where the hat notation indicates an approximation, $N$ is the number of degrees of freedom in the discretization, and $\{f_i(x)\}_{i=1}^N$ are basis functions satisfying $f_i(0)=f_i(1)=0$. The discretization allows measurements to be taken as a linear transformation of the degrees of freedom, which is given by $v_i(t) = \mathcal M_{ij} \theta_j(t)$, where $\mathcal M_{ij} = f_j(x_i)$ is the discretized measurement operator. To obtain a numerical solution, we compute the Galerkin weak form of Eq. \eqref{burgers_param} in space and plug in the discretization to obtain a system of ordinary differential equations. This reads:

\begin{equation*}
\begin{aligned}
    \sum_{j=1}^N\dot \theta_j \int_0^1 f_j f_i dx + \sum_{j=1}^N \sum_{k=1}^N \theta_j \theta_k \int_0^1 f_j \pd{f_k}{x} f_i dx + 10^{-\epsilon_1}\sum_{j=1}^N \theta_j \int_0^1 \pd{f_j}{x} \pd{f_i}{x} dx = \epsilon_2 \int_0^1 s(x,t) f_i dx,
\end{aligned}   
\end{equation*}

\noindent where we have used integration by parts to transfer derivatives from the discretized velocity to the test function in the diffusion term. Pre-computing and storing all matrix quantities defined by integrals of the basis functions, this expression is re-written as

\begin{equation*}
    M_{ij} \dot \theta_j + A_{ijk} \theta_j \theta_k + 10^{-\epsilon_1} K_{ij}\theta_j = \epsilon_2 F_i(t).
\end{equation*}

This is a system of ordinary differential equations whose solution can be approximated with time integration. Using the backward Euler method, the nonlinear system of equations governing the solution at the next time step is

\begin{equation}\label{backward_euler}
    \mathbf R( \boldsymbol \theta_{t+1} ; \boldsymbol \theta_t,\epsilon_1 , \epsilon_2) = \qty( \frac{\mathbf{M}}{\Delta t} + 10^{-\epsilon_1} \mathbf{K} ) \boldsymbol \theta_{t+1} + \mathbf{A} :( \boldsymbol \theta_{t+1} \otimes \boldsymbol \theta_{t+1}) - \frac{\mathbf{M}}{\Delta t} \boldsymbol \theta_t - \epsilon_2 \mathbf{F}_{t+1} = \mathbf{0},
\end{equation}

\noindent where $\Delta t$ is a fixed time step and---by starting with the initial condition $ \theta_{0,i}=\int_0^1u_0(x) f_i(x) dx$ and marching forward in time---the solution at the previous time step $\boldsymbol \theta_t$ is known. As a shorthand, we denote the tensor quantity $A_{ijk} \theta_j \theta_k$ with $\mathbf A : ( \boldsymbol \theta_{t+1} \otimes \boldsymbol \theta_{t+1})$. The total simulation time is $T$ and the number of time steps is $\mathcal P$, such that $\Delta t = T/\mathcal P$. This system of equations is solved with Newton's method, given by 

\begin{equation*}
    \boldsymbol \theta_{t+1}^{k+1} = \boldsymbol \theta_{t+1}^{k} - \qty( \pd{\mathbf R}{\boldsymbol  \theta}(\boldsymbol \theta_{t+1}^k))^{-1} \mathbf R(\boldsymbol \theta_{t+1}^{k}),
\end{equation*}

\noindent which continues until $\lVert \mathbf R (\boldsymbol \theta_{t+1}^{k})\rVert < \mathcal T$, where $\mathcal T$ is a user-defined convergence criterion. Even when the dependence of the residual on the solution at the previous time step and the model parameters $\boldsymbol \epsilon$ is omitted, it is understood to depend on these quantities.

\paragraph{} We first present the standard formulation of the inverse problem before proceeding to modify the governing equation to accommodate constraint forces. With the standard approach, the inverse problem is solved by minimizing the discretized time integral of the squared error of the solution with the data at the measurement positions:

\begin{equation}\label{inv_obj}
    \begin{aligned}
        \underset{\boldsymbol \epsilon, \boldsymbol \theta_1,\boldsymbol \theta_2,\dots} {\text{argmin }} \frac{\Delta t}{2} \sum_{i=1}^{\mathcal P} \lVert \mathcal M \boldsymbol \theta_i - \mathbf v_i \rVert^2 \\
        \text{s.t. }  \mathbf R( \boldsymbol \theta_{t+1} ; \boldsymbol \theta_t,\boldsymbol \epsilon) = \mathbf 0 \quad t=0,1,2,\dots, \mathcal P-1,
    \end{aligned}
\end{equation}

\noindent where $\mathbf v_i = \mathbf v(i\Delta t)$ is the measurement data evaluated on the time integration grid. We call the above objective function $z^{\text{INV}}$. We choose to treat the solution parameters as an explicit function of the model parameters through the constraint equation with sensitivity analysis. First, we assume that $\boldsymbol \theta(t) = \boldsymbol \theta(t;\boldsymbol \epsilon)$, meaning that the solution parameters are always computed through Eq. \eqref{backward_euler}. Second, noting this dependence, we write the gradient of the objective as

\begin{equation}\label{grad}
   \pd{}{\boldsymbol \epsilon} z^{\text{INV}}(\boldsymbol \epsilon) = \pd{}{\boldsymbol \epsilon} \frac{\Delta t}{2} \sum_{i=1}^{\mathcal P} \lVert \mathcal M \boldsymbol \theta_i - \mathbf v_i \rVert^2  = \Delta t\sum_{i=1}^{\mathcal P} (\mathcal M \boldsymbol \theta_i - \mathbf v_i ) \cdot \mathcal M \pd{\boldsymbol \theta_i}{\boldsymbol \epsilon}.
\end{equation}

To compute the sensitivity derivative $\partial \boldsymbol \theta_t / \partial \boldsymbol \epsilon $, we differentiate the discretized governing equation in continuous time with respect to the two model parameters:

\begin{equation*}
\begin{aligned} 
    M_{ij} \pd{}{t}\qty(\pd{\theta_j}{\epsilon_1}) + A_{ijk} \pd{\theta_j}{\epsilon_1} \theta_k + A_{ijk} \theta_j \pd{\theta_k}{\epsilon_1} - \ln( 10) 10^{-\epsilon_1}K_{ij}\theta_j + 10^{-\epsilon_1}K_{ij}\pd{\theta_j}{\epsilon_1} = 0,\\
     M_{ij} \pd{}{t}\qty(\pd{\theta_j}{\epsilon_2}) + A_{ijk} \pd{\theta_j}{\epsilon_2} \theta_k + A_{ijk} \theta_j \pd{\theta_k}{\epsilon_2} - F_i(t) = 0.  
\end{aligned}
\end{equation*}

A solution to these two systems of equations is approximated with backward Euler time integration on the same integration grid as the forward solution. As such, the governing system of equations for the sensitivity derivatives is 

\begin{equation}\label{discrete_sens}
    \begin{bmatrix}
        \frac{\mathbf M}{\Delta t} \pd{\boldsymbol \theta_{t+1}}{\epsilon_1} -\frac{\mathbf M}{\Delta t} \pd{\boldsymbol \theta_{t}}{\epsilon_1} +  \mathbf A : \pd{\boldsymbol \theta_{t+1}}{\epsilon_1} \otimes \boldsymbol \theta_{t+1} + \mathbf A :  \boldsymbol \theta_{t+1} \otimes \pd{\boldsymbol \theta_{t+1}}{\epsilon_1}  - \ln(10)10^{-\epsilon_1}\mathbf K \boldsymbol \theta_{t+1} + 10^{-\epsilon_1} \mathbf K \pd{\boldsymbol \theta_{t+1}}{\epsilon_1} \\
       \frac{\mathbf M}{\Delta t} \pd{\boldsymbol \theta_{t+1}}{\epsilon_2} -\frac{\mathbf M}{\Delta t} \pd{\boldsymbol \theta_{t}}{\epsilon_2} +  \mathbf A : \pd{\boldsymbol \theta_{t+1}}{\epsilon_2} \otimes \boldsymbol \theta_{t+1} + \mathbf A :  \boldsymbol \theta_{t+1} \otimes \pd{\boldsymbol \theta_{t+1}}{\epsilon_2} - \mathbf F_{t+1}
    \end{bmatrix} = \mathbf 0.
\end{equation}

Note that the initial conditions on the sensitivity derivatives is zero because the initial condition of the solution parameters does not change with the model parameters. This system is also solved at each time step with Newton's method. Thus, the gradient of the objective---which is now unconstrained---is computed with Eqs. \eqref{grad} and \eqref{discrete_sens} and supplied to the optimizer in search of a minimum error setting of the model parameters. This is the approach we take to find a solution to the standard formulation of the dynamic inverse problem given by Eq. \eqref{inv_obj}.

\paragraph{} To construct the ECFM approach to the inverse problem, we first introduce time-dependent constraint forces to the governing equation of Eq. \eqref{burgers_param}:

\begin{equation}\label{burgers_ecfm}
    \pd{w}{t} + w \pd{w}{x} = 10^{-\epsilon_1} \pdd{w}{x} + \epsilon_2s(x,t) + \sum_{i=1}^C \lambda_i(t) \Gamma(x-x_i), 
\end{equation}

\noindent where the boundary and initial conditions remain unchanged and $\Gamma(x-x_i)$ are constraint forces centered on the measurement locations. The spatial form of the constraint forces is chosen by the analyst. The constraint force magnitudes $\lambda_i(t)$ are additional unknowns used to enforce that the system agrees with the time-dependent measurement data. Computing the weak form of Eq. \eqref{burgers_ecfm}, we obtain 

\begin{equation*}
    M_{ij} \dot \theta_j + A_{ijk} \theta_j \theta_k + 10^{-\epsilon_1} K_{ij}\theta_j = \epsilon_2 F_i(t) + \Gamma_{ij}\lambda_j,
\end{equation*}

\noindent where $\Gamma_{ij} = \int_0^1 \Gamma(x-x_j) f_i dx$ is the force vector corresponding to the constraint forces. Again, using backward Euler time integration, the solution at the next time step is governed by 

\begin{equation}
    \qty( \frac{\mathbf{M}}{\Delta t} + 10^{-\epsilon_1} \mathbf{K}) \boldsymbol \theta_{t+1} + \mathbf{A} : \boldsymbol \theta_{t+1} \otimes \boldsymbol \theta_{t+1}  - \frac{\mathbf{M}}{\Delta t} \boldsymbol \theta_t - \epsilon_2 \mathbf{F}_{t+1} - \boldsymbol \Gamma \boldsymbol \lambda_{t+1} = \mathbf{0}.
\end{equation}

Recalling that the constraint forces enforce the data constraint, we solve an expanded system of equations to jointly determine the solution parameters and constraint force magnitudes:

\begin{equation}
    \mathbf{\tilde R}(\boldsymbol \theta_{t+1} , \boldsymbol \lambda_{t+1} ; \boldsymbol \theta_t ,  \epsilon_1 , \epsilon_2) = \begin{bmatrix}
        \qty( \frac{\mathbf{M}}{\Delta t} + 10^{-\epsilon_1} \mathbf{K} ) \boldsymbol \theta_{t+1} + \mathbf{A} : \boldsymbol \theta_{t+1} \otimes \boldsymbol \theta_{t+1} - \frac{\mathbf{M}}{\Delta t} \boldsymbol \theta_t - \epsilon_2 \mathbf{F}_{t+1} - \boldsymbol \Gamma \boldsymbol \lambda_{t+1} \\ \mathcal M \boldsymbol \theta_{t+1} - \mathbf{v}_{t+1}
    \end{bmatrix} = \begin{bmatrix}
        \mathbf{\tilde R}_1 \\ \mathbf{\tilde R}_2
    \end{bmatrix} = \mathbf{0}.
\end{equation}

In analogy to Eq. \eqref{ecfm}, the ECFM formulation of the time-dependent inverse problem is 

\begin{equation}\label{ecfm_obj}
    \begin{aligned}
        \underset{\boldsymbol \epsilon, \boldsymbol \theta_1,\boldsymbol \lambda_1 ,\boldsymbol \theta_2,\boldsymbol \lambda_2,\dots} {\text{argmin }} \frac{\Delta t}{2} \sum_{i=1}^{\mathcal P} \lVert \boldsymbol \lambda_i\rVert^2 \\
        \text{s.t. }  \mathbf {\tilde R}( \boldsymbol \theta_{t+1} ,\boldsymbol \lambda_{t+1}; \boldsymbol \theta_t,\boldsymbol \epsilon) = \mathbf 0 \quad t=0,1,2,\dots, \mathcal P-1.
    \end{aligned}
\end{equation}

We call the objective function measuring the constraint force magnitude $z^{\text{ECFM}}$. We use sensitivity analysis to rid the optimization problem of the constraint, though this is more delicate in the case of ECFM. We treat the constraint forces as explicit functions of the model parameters $\boldsymbol \epsilon$, and eliminate the dependence of the optimization problem on the solution parameters $\boldsymbol \theta$. To accomplish this, we first solve for the constraint forces as a function of the model parameters using Newton's method at each time step:

\begin{equation}\label{ecfm_solve}
    \begin{bmatrix}
        \boldsymbol \theta_{t+1 } \\ \boldsymbol \lambda_{t+1}
    \end{bmatrix}^{k+1} = \begin{bmatrix}
        \boldsymbol \theta_{t+1 } \\ \boldsymbol \lambda_{t+1}
    \end{bmatrix}^{k} - \begin{bmatrix}
        \partial \mathbf{\tilde R}_1 / \partial \boldsymbol \theta & \partial \mathbf{\tilde R}_1  / \partial \boldsymbol \lambda \\ \partial \mathbf{\tilde R}_2 / \partial \boldsymbol \theta & \partial \mathbf{\tilde R}_2 / \partial \boldsymbol \lambda
    \end{bmatrix}^{-1} \begin{bmatrix}
        \mathbf{\tilde R}_1(\boldsymbol \theta_t^k, \boldsymbol \lambda^k_t) \\ \mathbf{\tilde R}_2(\boldsymbol \theta^k_t, \boldsymbol \lambda^k_t)
    \end{bmatrix},
\end{equation}

\noindent which continues until $\lVert \mathbf{\tilde R}(\boldsymbol \theta_t, \boldsymbol \lambda_t) \rVert < \mathcal T$, where $\mathcal T$ is a user-defined convergence criterion. Though we do not write it explicitly, it is understood that the residual system depends on the model parameters $\boldsymbol \epsilon$ and that the Jacobian matrix is evaluated at the current Newton iterate of the solution parameters and constraint forces.  Eq. \eqref{ecfm_solve} enforces that the constraint forces are explicit functions of the model parameters. Noting this dependence, we write the gradient of the objective in Eq. \eqref{ecfm_obj} as 

\begin{equation}\label{ecfm_grad}
    \pd{}{\boldsymbol \epsilon} z^{\text{ECFM}}(\boldsymbol \epsilon) = \pd{}{\boldsymbol \epsilon} \frac{\Delta t}{2} \sum_{i=1}^{\mathcal P} \lVert \boldsymbol \lambda_i\rVert^2 = \Delta t \sum_{i=1}^{\mathcal P} \boldsymbol \lambda_i \cdot \pd{\boldsymbol \lambda_i}{\boldsymbol \epsilon}.
\end{equation}

In order to compute the sensitivity of the constraint force, we differentiate the discretized governing equation in continuous time:

\begin{equation}\label{sens_ecfm}
\begin{aligned} 
    M_{ij} \pd{}{t}\qty(\pd{\theta_j}{\epsilon_1}) + A_{ijk} \pd{\theta_j}{\epsilon_1} \theta_k + A_{ijk} \theta_j \pd{\theta_k}{\epsilon_1} - \ln( 10) 10^{-\epsilon_1}K_{ij}\theta_j + 10^{-\epsilon_1}K_{ij}\pd{\theta_j}{\epsilon_1} - \Gamma_{ij} \pd{\lambda_j}{\epsilon_1}= 0,\\
     M_{ij} \pd{}{t}\qty(\pd{\theta_j}{\epsilon_2}) + A_{ijk} \pd{\theta_j}{\epsilon_2} \theta_k + A_{ijk} \theta_j \pd{\theta_k}{\epsilon_2} - F_i(t) - \Gamma_{ij} \pd{\lambda_j}{\epsilon_2} = 0.  
\end{aligned}
\end{equation}

Now, we must also differentiate the continuous time constraint equation. This introduces two additional contributions to the system of equations governing the sensitivity of the solution parameters and constraint forces:

\begin{equation*}
    \mathcal M \pd{\boldsymbol \theta}{\epsilon_1} = \mathcal M \pd{\boldsymbol \theta}{\epsilon_2} = \mathbf 0.
\end{equation*}

Using the results of the forward solve of Eq. \eqref{ecfm_solve}, approximating the time derivatives in Eq. \eqref{sens_ecfm} with the backward Euler method, and augmenting the resulting system with the constraint equations, we obtain a scheme to solve for the constraint force sensitivities at each time step. This is given by

\begin{equation}\label{discrete_sens_ecfm}
    \begin{bmatrix}
        \frac{\mathbf M}{\Delta t} \pd{\boldsymbol \theta_{t+1}}{\epsilon_1} -\frac{\mathbf M}{\Delta t} \pd{\boldsymbol \theta_{t}}{\epsilon_1} +  \mathbf A : \pd{\boldsymbol \theta_{t+1}}{\epsilon_1} \otimes \boldsymbol \theta_{t+1} + \mathbf A :  \boldsymbol \theta_{t+1} \otimes \pd{\boldsymbol \theta_{t+1}}{\epsilon_1}  - \ln(10)10^{-\epsilon_1}\mathbf K \boldsymbol \theta_{t+1} + 10^{-\epsilon_1} \mathbf K \pd{\boldsymbol \theta_{t+1}}{\epsilon_1} - \boldsymbol \Gamma \pd{\boldsymbol \lambda_{t+1}} {\epsilon_1}\\
        \mathcal M \pd{\boldsymbol \theta_{t+1}}{\epsilon_1}\\
       \frac{\mathbf M}{\Delta t} \pd{\boldsymbol \theta_{t+1}}{\epsilon_2} -\frac{\mathbf M}{\Delta t} \pd{\boldsymbol \theta_{t}}{\epsilon_2} +  \mathbf A : \pd{\boldsymbol \theta_{t+1}}{\epsilon_2} \otimes \boldsymbol \theta_{t+1} + \mathbf A :  \boldsymbol \theta_{t+1} \otimes \pd{\boldsymbol \theta_{t+1}}{\epsilon_2} - \mathbf F_{t+1} - \boldsymbol \Gamma \pd{\boldsymbol \lambda_{t+1}}{\epsilon_2} \\
        \mathcal M \pd{\boldsymbol \theta_{t+1}}{\epsilon_2}
    \end{bmatrix} = \mathbf 0
\end{equation}

Note that the initial conditions on both the sensitivity derivatives of the solution parameters and the constraint forces are zero because the initial velocity does not change with the model parameters, and thus no constraint force is required to enforce agreement with the initial measurements. When the constraint force objective is computed at a given $\boldsymbol \epsilon$ through Eq. \eqref{ecfm_solve}, and its gradient is computed with Eqs. \eqref{ecfm_grad} and \eqref{discrete_sens_ecfm}, the minimization of the constraint force magnitude is without constraints. We now compare the performance of the constraint force formulation of the inverse problem to the standard approach.

\begin{figure}[h]
\centering
\includegraphics[width=0.99\textwidth]{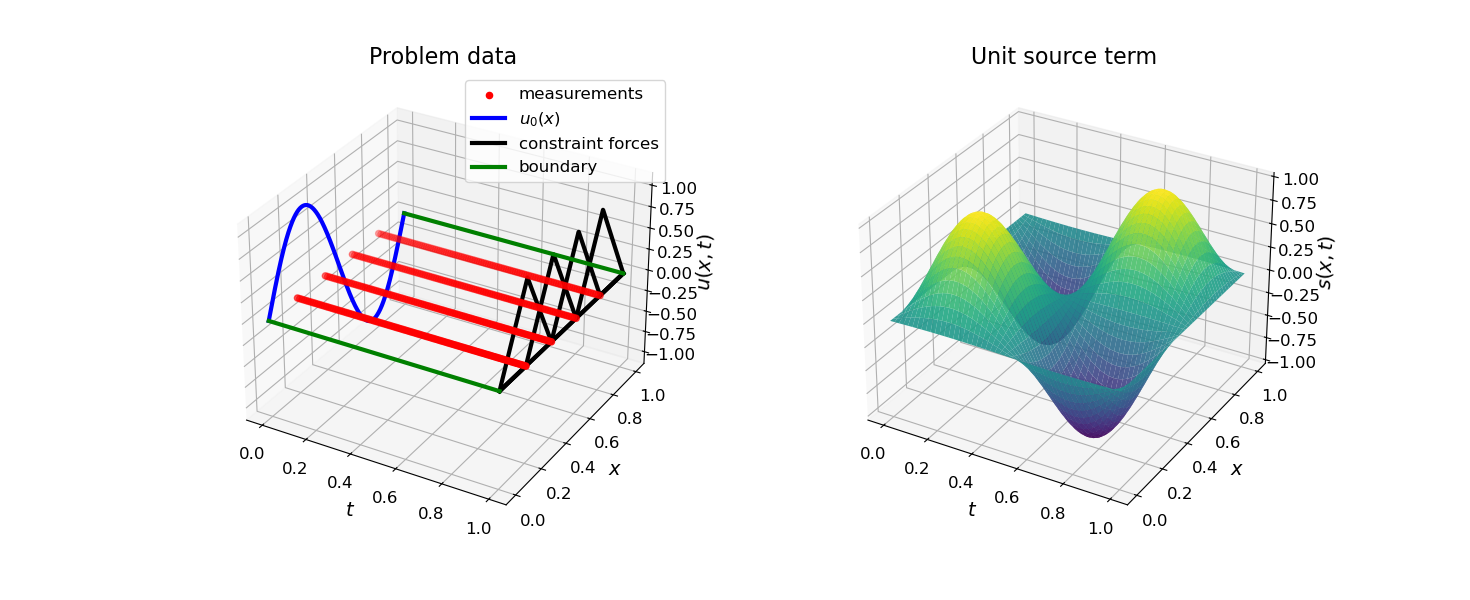}
\caption{Problem setup for extracting the viscosity and source term magnitude from flow data generated by solving Burgers' equation.}
\label{setup}
\end{figure}

\begin{figure}[h]
\centering
\includegraphics[width=0.99\textwidth]{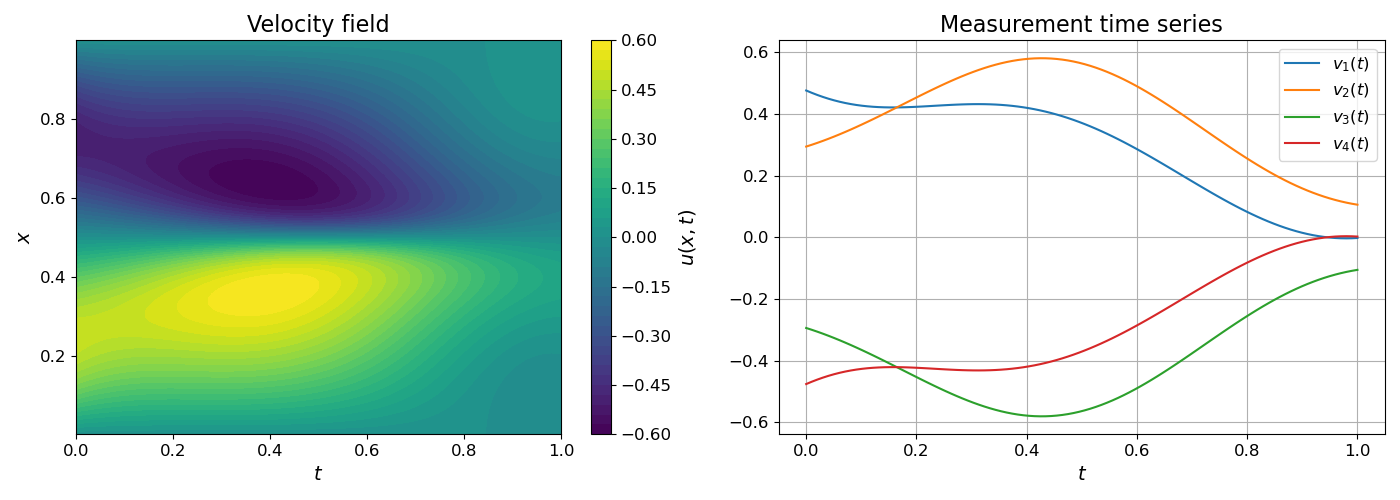}
\caption{True solution and the time series of the $C=4$ measurements taken, which are used in the inverse problem to recover the viscosity and source term magnitude.}
\label{sol}
\end{figure}

\paragraph{} In this example, we choose to work with zero measurement noise and an inverse problem in which the parameterized model is consistent with the measurement data. Thus, we expect a successful method to exactly recover the parameters of the true model if sufficient data is provided. We take the initial condition to be $u_0(x)=\sin(2\pi x)$, which gives rise to the well-known shock pattern when $\nu=0$. The solution is discretized with basis functions given by $f_i(x) = \sin(i \pi x)$ for $i=1,\dots,N$, which ensures that the homogeneous Dirichlet boundary conditions are satisfied automatically. We take $N=50$ shape functions in the discretization, $C=4$ evenly spaced measurement points, simulate the problem over a time interval $T=2$, using $\mathcal P=100$ time integration points. Eqs. \eqref{burgers} and \eqref{burgers_param} show that the true magnitude of the source term is $\epsilon_2=1$. We take the space-time form of the source term to be $s(x,t)=\sin(2\pi x) \sin( 2\pi t)$. The true viscosity is $\nu=1.78 \times 10^{-2}$, which corresponds to a viscosity model parameter of $\epsilon_1=1.75$. The constraint forces are one-dimensional finite element hat functions centered on the measurement positions. See Figure \ref{setup} to visualize the initial condition, the measurement/time integration grid, the constraint forces, and the source term. Figure \ref{sol} shows the true velocity field and the noiseless measurement data taken from it. The convergence criterion for the Newton solve in both the standard inverse problem and the constraint force formulation is $\mathcal T = 10^{-6}$. For both the standard inverse problem and the ECFM formulation, we use ADAM optimization with a learning rate of $10^{-2}$ to find a minimum of the unconstrained objective with gradients manually computed using sensitivity analysis. For both problems, ADAM optimization is run for $250$ epochs.

\begin{figure}[h]
\centering
\includegraphics[width=0.99\textwidth]{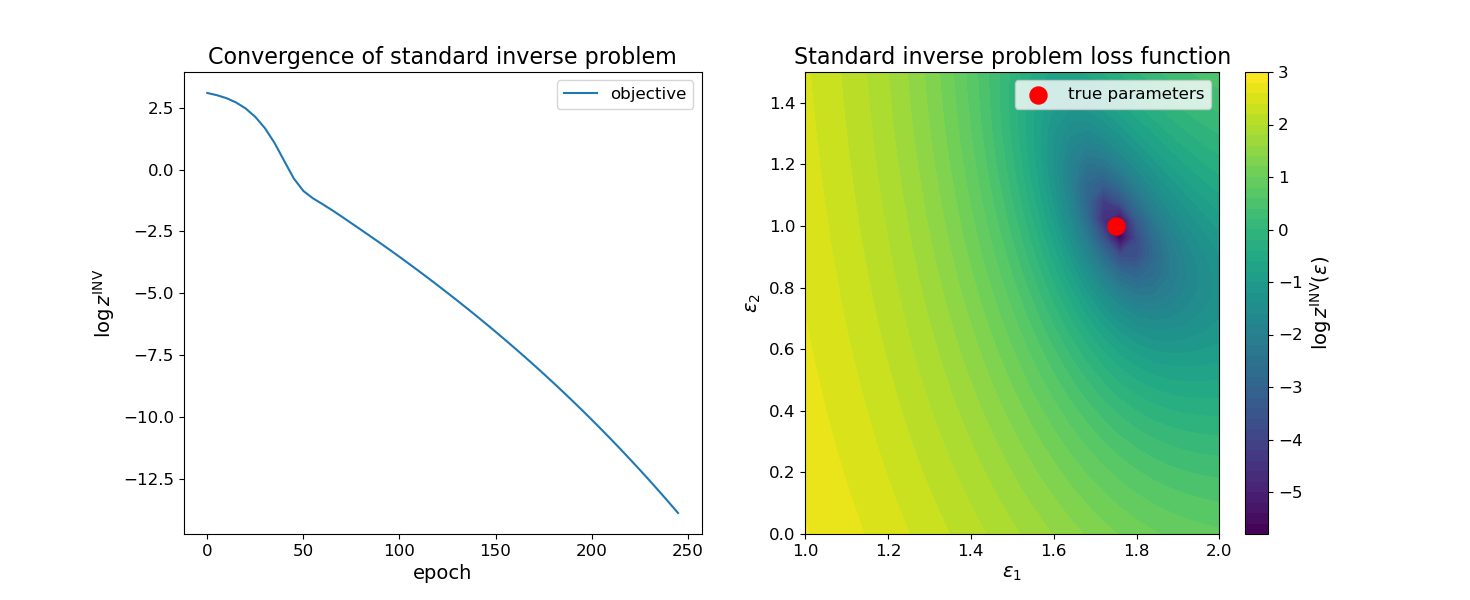}
\caption{The standard inverse problem converges to a zero-error solution because the parameterized model is consistent with the measurement data (left). The loss surface is convex in the vicinity of the true solution, indicating an identifiable inverse problem (right).}
\label{inv}
\end{figure}

\begin{figure}[h]
\centering
\includegraphics[width=0.99\textwidth]{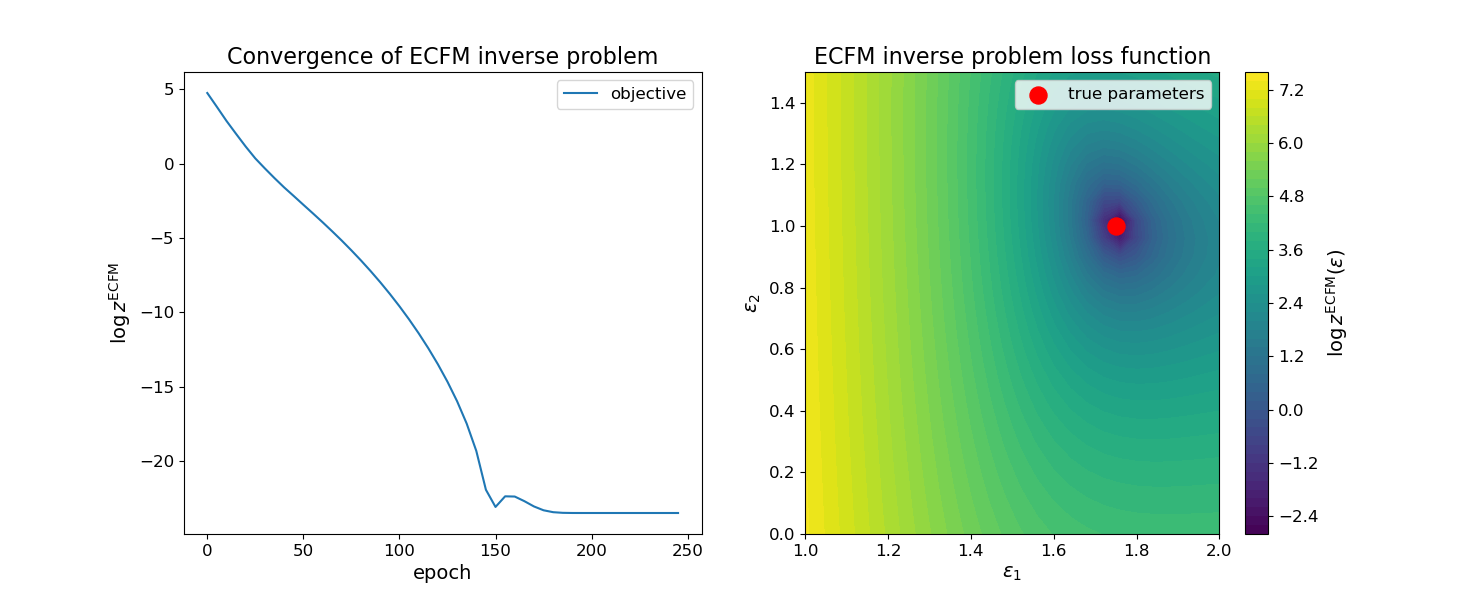}
\caption{The constraint force formulation of the inverse problem converges to model parameters with zero associated constraint force because the parameterized model is consistent with the measurement data (left). The loss surface representing the constraint force magnitude is also convex in the vicinity of the true solution, indicating an identifiable inverse problem (right).}
\label{ecfm_fig}
\end{figure}

\paragraph{} Figure \ref{inv} shows the results of the standard inverse problem as well as the corresponding loss surface. After the allotted $250$ epochs, the recovered parameters are $\epsilon_1^{\text{INV}}=1.750$ and $\epsilon_2^{\text{INV}}=1.000$, exactly matching the true values to three decimal places. The loss surface is convex in the vicinity of the solution, indicating a well-posed inverse problem. Compare this to Figure \ref{ecfm_fig}, where the convergence of the constraint force objective and the corresponding loss surface are shown. The recovered parameters with ECFM are also $\epsilon_1^{\text{ECFM}}=1.750$ and $\epsilon_2^{\text{ECFM}}=1.000$. We remark that the loss surface arising from ECFM is better-conditioned than that of the standard inverse problem. Computing the Hessian matrix of the objective with finite differences around the solution, we obtain a condition number of $8.11$ for the standard inverse problem, and $3.46$ for the ECFM formulation. While this is interesting to note, we offer no opinion as to whether ECFM inverse problems are better-conditioned in general. We expect the standard and ECFM approaches to provide equivalent estimates of the model parameters when the parameterized model is consistent with the data. This is because a zero error solution in the sense of displacement/velocity is also a zero solution in the sense of constraint forces.

\subsection{Noisy measurements}

\paragraph{} With this example, we introduce a new approach to handling measurement noise with ECFM. In \cite{rowan_physics-informed_2025}, noisy measurements are treated with inequality constraints whose range represents a confidence interval based on the known distribution of the noise. This method for enforcing inequality constraints draws inspiration from the slack variable approach familiar from Lagrange multipliers \cite{ruszczynski_nonlinear_2006}. But, because the form of the constraint forces is controlled by the analyst, the governing system of equations no longer follows from stationarity of a Lagrange function. Thus, with this approach, a standard inequality-constrained optimizer cannot be recycled to treat inequality constraints with ECFM. Instead, the analyst must develop the routine to handle slack variables and the analogue of the ``dual feasibility'' requirement on the constraint forces from scratch \cite{rowan_physics-informed_2025}. We take this to be a prohibitively complex and/or intrusive strategy to deal with measurement noise. Yet, when using constraint forces in the presence of measurement noise, it is important to handle the noise in a principled way. Unlike the squared error objective of Eq. \eqref{traditional}, which is robust to measurement noise, the constraint force objective of Eq. \eqref{ecfm} should not be used with noise, given that equality constraints on the measurement data are enforced. Depending on the measurement positions, very small deviations of the solution from the measurement data may require large constraint forces to correct. For example, if a noisy measurement is taken near a Dirichlet boundary of the system, large constraint forces are required in order for the solution to track the measurement noise, as the system is stiff near the boundary. Given that the magnitude of the constraint force provides insight into the consistency of the parameterized model with the measurement data, such a situation would give a false impression of inconsistency, whereas in reality, the constraint forces would be small or zero if noise was accounted for. Thus, we must devise a new approach to handling noise which is more tractable than that of \cite{rowan_physics-informed_2025}.


\paragraph{} Before we begin the discussion of measurement noise, we first introduce the physics of the test problem. We work with a scalar-valued nonlinear static PDE called the Fisher-Kolmogorov-Petrovsky-Psikunov (Fisher-KPP) equation, which is a nonlinear reaction-diffusion PDE used to model population growth and traveling wave fronts \cite{rohrhofer_approximating_2025}. Steady-states of the Fisher-KPP equation are governed by

\begin{equation}\label{kpp}
\begin{aligned}
    \mu\nabla^2 u(\mathbf{x}) + r u(\mathbf{x})(1 - u(\mathbf{x})) + s(\mathbf{x}) = 0, \quad \mathbf{x} \in \Omega, \\
    u(\mathbf{x}) = 0, \quad \mathbf{x} \in \partial \Omega,
\end{aligned}
\end{equation}

\noindent where $\mu$ is the diffusion rate and $r$ is the reaction rate. We assume homogeneous Dirichlet boundaries, and we set $\mu=1/2$ and $r=1$ going forward. The nature of the measurements of the true system will be discussed shortly. The parameterized boundary value problem used in the inverse problem is 

\begin{equation}\label{kpp_param}
\begin{aligned}
    \frac{1}{2}\nabla^2 w(\mathbf{x}) +  w(\mathbf{x})(1 - w(\mathbf{x})) + b(\mathbf{x};\boldsymbol \epsilon) = 0, \quad \mathbf{x} \in \Omega, \\
    w(\mathbf{x}) = 0, \quad \mathbf{x} \in \partial \Omega.
\end{aligned}
\end{equation}

We do not assume that there is a parameter setting such that $b(\mathbf x;\boldsymbol \epsilon) = s(\mathbf x)$. In other words, we allow that the parameterized model is inconsistent with the true model, and thus, potentially inconsistent with the measurement data. The solution is discretized with $\hat w(\mathbf x; \boldsymbol \theta) = \sum_{i=1}^N \theta_i f_i (\mathbf x)$, and the source term is discretized with the same basis functions as $b(\mathbf x;\boldsymbol \epsilon) = \sum_{\alpha=1}^M \epsilon_{\alpha} f_{\alpha}(\mathbf x)$ where $N\neq M$ in general. For clarity, we use Greek indices to distinguish model parameters from the solution parameters with Latin indices. Like the previous example, we assume that the basis functions respect the homogeneous Dirichlet boundaries by construction. 

\paragraph{} First, we formulate the standard inverse problem, which requires no modifications to be used in the setting of noisy measurements. We then proceed to introduce a novel technique for handling noise with ECFM. Plugging in the discretization of the solution and the source term and computing the Galerkin weak form of Eq. \eqref{kpp_param}, we obtain

\begin{equation*}
    \frac{1}{2}\qty( \int_{\Omega} \nabla f_i \cdot \nabla f_j d\Omega) \theta_j - \qty(\int_{\Omega}f_i f_j d\Omega) \theta_j + \qty(\int_{\Omega} f_i f_j f_k d\Omega) \theta_j \theta_k - \qty(\int_{\Omega} f_i f_{\alpha}d\Omega) \epsilon_\alpha = 0.
\end{equation*}

Computing all integrals of the basis functions, this expression can be written in symbolic form as

\begin{equation*}
    \mathbf R(\boldsymbol \theta ; \boldsymbol \epsilon) = (\mathbf K - \mathbf M)\boldsymbol \theta + \mathbf A : \boldsymbol \theta \otimes \boldsymbol \theta - \mathbf F \boldsymbol \epsilon =\mathbf 0,
\end{equation*}




\noindent where the notation $\mathbf A : \boldsymbol \theta \otimes \boldsymbol \theta$ denotes contraction on the last two indices, i.e., $A_{ijk} \theta_j \theta_k$. With the discretized measurement operator $\mathcal M$ and measurement data $\mathbf v$, the standard formulation of the inverse problem is

\begin{equation}\label{kpp_inv}
    \begin{aligned}
        \underset{\boldsymbol \epsilon, \boldsymbol \theta}{\text{argmin }} \frac{1}{2} \lVert \mathcal M \boldsymbol \theta - \mathbf v\rVert^2 \\
        \text{s.t. }\mathbf R(\boldsymbol \theta ; \boldsymbol \epsilon) = \mathbf 0.
    \end{aligned}
\end{equation}

\paragraph{} Eq. \eqref{kpp_inv} is the standard formulation of the static inverse problem. Now, ECFM requires the introduction of constraint forces into the governing equation of Eq. \eqref{kpp_param}. The constraint force source term corresponding to $C$ measurements is $\sum_{i=1}^C \lambda_i \Gamma(\mathbf x - \mathbf x_i)$. Introducing the constraint forces, plugging in the discretization of the state, and computing the weak form, the discretized ECFM system is 

\begin{equation*}
    \mathbf{\tilde R}(\boldsymbol \lambda , \boldsymbol \theta ; \boldsymbol \epsilon) = (\mathbf K - \mathbf M) \boldsymbol \theta + \mathbf A : \boldsymbol \theta \otimes \boldsymbol \theta - \mathbf F \boldsymbol \epsilon - \boldsymbol \Gamma \boldsymbol \lambda = \mathbf 0,
\end{equation*}

\noindent where the columns of $\boldsymbol \Gamma$ are the force vectors corresponding to the constraint forces at each of the measurement positions. In this example, we take the constraint forces to be radial basis functions centered at measurement points, i.e., $\Gamma(\mathbf x - \mathbf x_i) = \frac{\omega}{\pi}\exp(-\omega\lVert  \mathbf x - \mathbf x_i\rVert^2)$, where $\omega$ controls the width of the radial basis function. The constraint force vector is thus given by 

\begin{equation*}
    \Gamma_{ij} = \int_{\Omega}  f_i(\mathbf x) \Gamma( \mathbf x - \mathbf x_j)d\Omega.
\end{equation*}

Referring to the terminology of \cite{rowan_physics-informed_2025}, the question we face is: how do we define the analogue of the ``inner loop'' problem, which determines the constraint forces as a function of the model parameters $\boldsymbol \epsilon$? As discussed above, if the constraint forces are determined by the equality constraint $\mathcal M \boldsymbol \theta - \mathbf v = \mathbf 0$, spuriously large constraint forces may be required for the system to track the measurement noise. We propose the following strategy to avoid this: find the minimum constraint force such that the sample mean and variance of the discrepancy $\mathcal M \boldsymbol \theta - \mathbf v$ fall within a $(1-\alpha)$\% confidence interval around their true values. This ensures that the discrepancy replicates the known statistical properties of the measurement operation. The constraint forces thus correct model misspecification---as the statistics of the discrepancy of a misspecified model will not resemble noise in general---but they do not needlessly force the system to follow the noisy data. Enforcing that the sample mean and variance of the discrepancy fall within a confidence interval gives rise to an inequality-constrained optimization problem:

\begin{equation*}
\begin{aligned}
\underset{\boldsymbol \epsilon, \boldsymbol \lambda, \boldsymbol \theta}{\text{argmin}}
\quad & \frac{1}{2}\lVert \boldsymbol \lambda \rVert^{2} \\[6pt]
\text{s.t.}\quad 
& \tilde{\mathbf R}(\boldsymbol \lambda, \boldsymbol \theta;\boldsymbol \epsilon) = \mathbf 0, \\[6pt]
& \ell_1 \;\le\; 
\frac{1}{C}\sum_{i=1}^C ( \mathcal M_{ij}\theta_j - v_i )
\;\le\; \ell_2, \\[6pt]
& p_1 \;\le\;
\frac{1}{C-1}
\left(
 (\mathcal M_{ij}\theta_j - v_i)( \mathcal M_{ik}\theta_k - v_i)
 - \frac{1}{C}
   \left( \sum_{i=1}^C \mathcal M_{ij}\theta_j - v_i \right)^2
\right)
\;\le\; p_2 ,
\end{aligned}
\end{equation*}

\noindent where $\ell_1$, $\ell_2$, $p_1$, and $p_2$ are bounds based on the sampling distribution of the mean and variance which we discuss shortly. The equality constraint enforces that the system is in equilibrium, the first inequality constraint enforces that the sample mean of the discrepancy is in a given range, and the second inequality constraint enforces that the sample variance is in a given range. We write explicit summations when indices are not repeated, but rely on the summation convention for repeated indices. Calling the discrepancy $e_i = \mathcal M_{ij} \theta_j - v_i$, this equation can be written in a standard form as:

\begin{equation}\label{ecfm_noise}
\begin{aligned}
\underset{\boldsymbol \epsilon, \boldsymbol \lambda, \boldsymbol \theta}{\text{argmin}}
\quad & \frac{1}{2}\lVert \boldsymbol \lambda \rVert^{2} \\[6pt]
\text{s.t.}\quad 
& \tilde{\mathbf R}(\boldsymbol \lambda, \boldsymbol \theta;\boldsymbol \epsilon) = \mathbf 0, \\[6pt]
& \frac{1}{C} \sum_{i=1}^C e_i - \ell_1 \;\geq\ 0, \\[6pt]
& \ell_2 - \frac{1}{C}\sum_{i=1}^C e_i \;\geq\ 0, \\[6pt]
& \frac{1}{C-1}\qty( \sum_{i=1}^C e_i^2 - \frac{1}{C}\qty(\sum_{i=1}^C e_i)^2) - p_1 \;\geq\ 0 ,\\[6pt]
& p_2 - \frac{1}{C-1}\qty( \sum_{i=1}^C e_i^2 - \frac{1}{C}\qty(\sum_{i=1}^C e_i)^2) \;\geq\ 0.
\end{aligned}
\end{equation}

We remark that Eq. \eqref{ecfm_noise} breaks the structure of an ``inner'' and ``outer'' loop problem with which ECFM was originally presented. This hierarchical structure was the consequence of viewing the solution parameters and constraint forces as explicit functions of the model parameters, but this was also the cause of the customized procedure for inequality constraints. If we forego the inner-outer loop structure, we are free to use a standard inequality-constrained optimizer to solve Eq. \eqref{ecfm_noise}. To expedite the optimization process, gradients of the objective and constraints are computed analytically. Because the three sets of parameters have no explicit dependence on each other, it is not necessary to perform sensitivity analysis when computing gradients. 

\paragraph{} We now discuss the determination of the bounds for the inequality constraints on the sample mean and variance. Per Section 2, we assume that the noise is normally distributed with zero mean and known variance $\sigma^2$. A result from elementary statistics \cite{casella_statistical_2024} states that the ``sampling distribution'' of the sample mean is given by 

\begin{equation*}
    \frac{1}{C} \sum_{i=1}^C e_i \sim N\qty(0,\frac{\sigma^2}{C}),
\end{equation*}

\noindent which corresponds to a standard deviation of $\sigma/ \sqrt{C}$. We thus determine a $(1-\alpha)$\% confidence interval by finding the limits $\ell_1$ and $\ell_2$ such that $P\qty(\ell_1 < \frac{1}{C} \sum_{i=1}^C e_i < \ell_2) = 1-\alpha$. For normal distributions, this is accomplished with the help of the standard deviation and the z-score. This gives rise to:

\begin{equation*}
    \ell_1 = z_{\alpha/2} \frac{\sigma}{\sqrt{C}}, \quad \ell_2=-\ell_1,
\end{equation*}

\noindent where the z-score is defined implicitly through $\alpha/2=\int_{-\infty}^{z_{\alpha/2}} \exp(-x^2/2) /\sqrt{2\pi }dx$. The case of the sampling distribution of the sample variance is a bit more complicated. Another result from elementary statistics states that 

\begin{equation*}
    \frac{1}{C-1}\qty( \sum_{i=1}^C e_i^2 - \frac{1}{C}\qty(\sum_{i=1}^C e_i)^2) \sim \frac{\sigma^2}{C-1} \chi^2(C-1),
\end{equation*}

\noindent where $\chi^2(C-1)$ is a chi-squared distribution with $C-1$ degrees of freedom. We use the analogue of z-scores for chi-squared distribution, which are given by $\chi_{\alpha/2,C-1}$ and $\chi_{1-\alpha/2,C-1}$ respectively. These quantities are inverses of the cumulative distribution of the chi-squared random variable evaluated at $\alpha/2$ and $1-\alpha/2$. The limits of the confidence interval for the sample variance are thus

\begin{equation*}
    p_1 = \chi_{\alpha/2,C-1}\frac{\sigma^2}{C-1}, \quad p_2 =\chi_{1-\alpha/2,C-1}\frac{\sigma^2}{C-1}.
\end{equation*}

To clarify, the sampling distributions of the sample mean and sample variance of the discrepancy only follow these relations when they arise purely from measurement noise, and not model misspecification. Hence, the constraints in Eq. \eqref{ecfm_noise} are used to enforce the validity of the assumptions used in constructing these confidence intervals.

\paragraph{} We now perform a comparison of the standard inverse problem and the ECFM formulation on the Fisher-KPP equation. We take the domain to be the unit square $\Omega=[0,1]^2$ and the basis functions to be $N=100$ tensor products of univariate sines $\sin( i \pi x_1)$ and $\sin( j \pi x_2)$. The source term is discretized with $M=16$ basis functions of the same form. We take $C=225$ uniformly spaced measurements on the interior of the unit square domain and $\alpha=0.05$, giving rise to a $95$\% confidence interval on the sample mean and variance. This corresponds to $z_{0.025} = -1.96$, $\chi_{0.025,224}=184.44$, and $\chi_{0.975,224}=267.35$. The true source term is given by 

\begin{equation*}
    s(\mathbf x) = \begin{cases}
        100, \quad 0.25 \leq x_1,x_2 \leq 0.75,\\
        0, \quad \text{otherwise}.
    \end{cases}
\end{equation*}

\begin{figure}[h]
\centering
\includegraphics[width=0.99\textwidth]{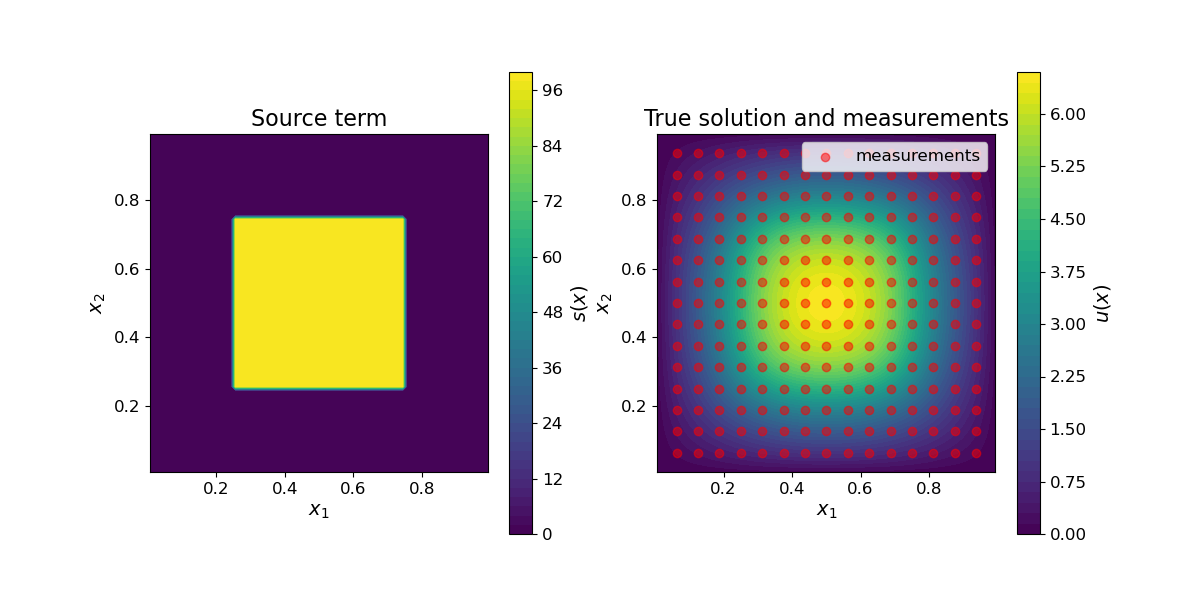}
\caption{The parameterized model is inconsistent with the measurement data because the true source term cannot be recovered by a finite sine series (left). The true solution is generated by approximating a solution to Eq. \eqref{kpp} with the discontinuous source term (right). Measurement data is taken at $C=225$ locations with normally distributed noise of a given mean and variance.}
\label{measure}
\end{figure}

\begin{figure}[h]
\centering
\includegraphics[width=0.99\textwidth]{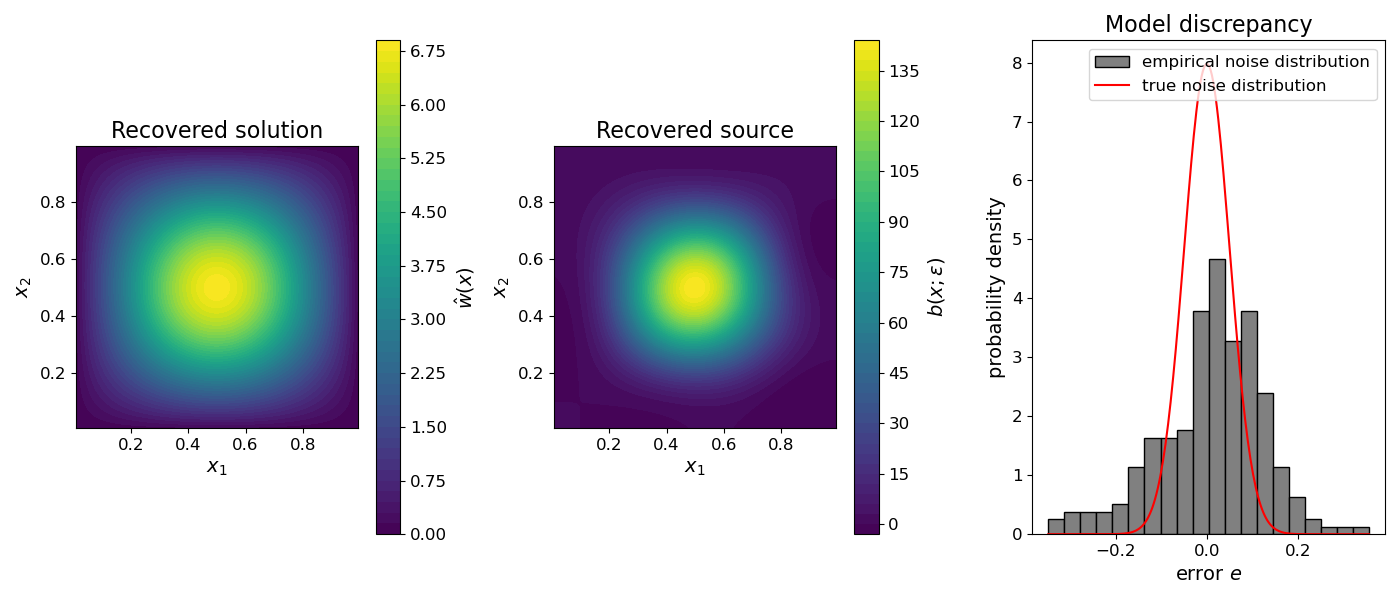}
\caption{The recovered solution obtains a small relative error with the true solution (left) despite the recovery of the source term being inaccurate (center). Though the relative error of the solution is small, the empirical distribution of the error does not reflect the known distribution of measurement noise (right).}
\label{kpp_inv_fig}
\end{figure}

The parameterized model of Eq. \eqref{kpp_param} is thus misspecified, as a finite combination of sines cannot exactly represent the discontinuous source term. The standard deviation of the measurement noise is given by $\sigma=5 \times 10^{-2}$. The measurement data is generated by approximating a solution to Eq. \eqref{kpp} with the same discretization and the true source term. See Figure \ref{measure} for the source term, the corresponding solution field, and the measurement grid. For both objectives, we use the scipy ``minimize'' optimizer with sequential least squares quadratic programming (SLSQP) to solve the constrained optimization problem. SLSQP iteratively approaches a stationary point of the Lagrange function using a quadratic approximation of the objective with the BFGS Hessian, a linear approximation of the constraints, and line search. We provide analytical gradients of both the objective and the constraints. Figure \ref{kpp_inv_fig} shows the results of the standard inverse problem. The relative error between the true and recovered solution is 

\begin{equation*}
    \mathcal E^{\text{INV}} = \frac{\int_{\Omega} | u(\mathbf x) - w^{\text{INV}}(\mathbf x) | d\Omega}{\int_{\Omega}  | u(\mathbf x) | d\Omega} = 3.2 \times 10^{-2}.
\end{equation*}

The recovered solution obtains low relative error despite the source term being inaccurately recovered. This is a consequence of the Laplacian operator damping the high-frequency basis functions which contribute to sharpening the boundaries of the square region of the source term. However, the statistical characteristics of the errors $\mathbf e = \mathcal M \boldsymbol \theta - \mathbf v$ do not match the known distribution of noise. At a $95$\% confidence level, the null hypothesis that the errors were drawn from a distribution with variance $\sigma^2$ is rejected. This disparity between the statistics of the error and the statistics of the noise is a consequence of model misspecification.

\begin{figure}[h]
\centering
\includegraphics[width=0.99\textwidth]{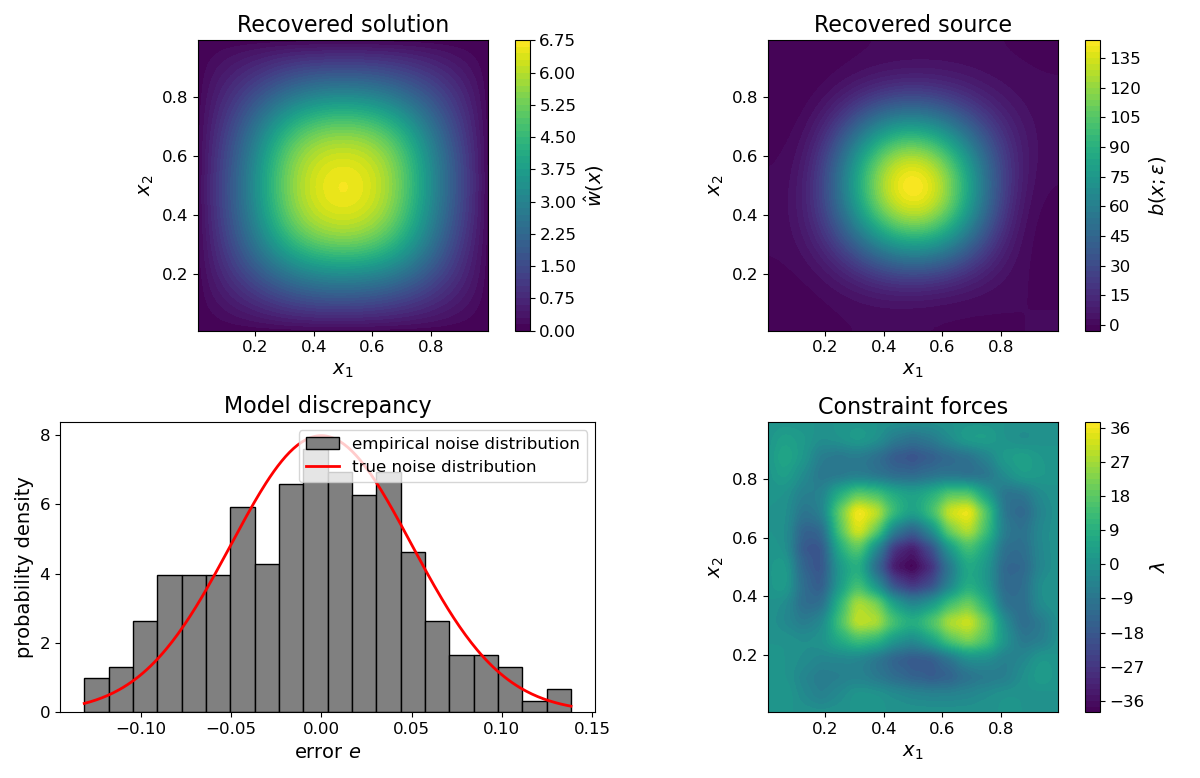}
\caption{The ECFM approach to inverse problems with known measurement noise. Constraint forces are introduced into the system to enforce that the discrepancy matches the statistical properties of the noise. The constraint forces are non-zero and exhibit clear structure, indicating missing physics in the parameterized model.}
\label{ecfm_results}
\end{figure}

\paragraph{} Using the bounds $\ell_1$, $\ell_2$, $p_1$, and $p_2$ obtained from the sampling distributions of the sample mean and variance, we use the ECFM formulation of the inverse problem to recover the model parameters. The results of solving the optimization problem of Eq. \eqref{ecfm_noise} are shown in Figure \ref{ecfm_results}. By construction, the first two moments of the discrepancy fall within their $95$\% confidence intervals. The relative error of the recovered solution with the true solution is given by

\begin{equation*}
    \mathcal E^{\text{ECFM}} = \frac{\int_{\Omega} | u(\mathbf x) - w^{\text{ECFM}}(\mathbf x) | d\Omega}{\int_{\Omega}  | u(\mathbf x) |d\Omega} = 1.1 \times 10^{-2}.
\end{equation*}

\begin{figure}[h]
\centering
\includegraphics[width=0.95\textwidth]{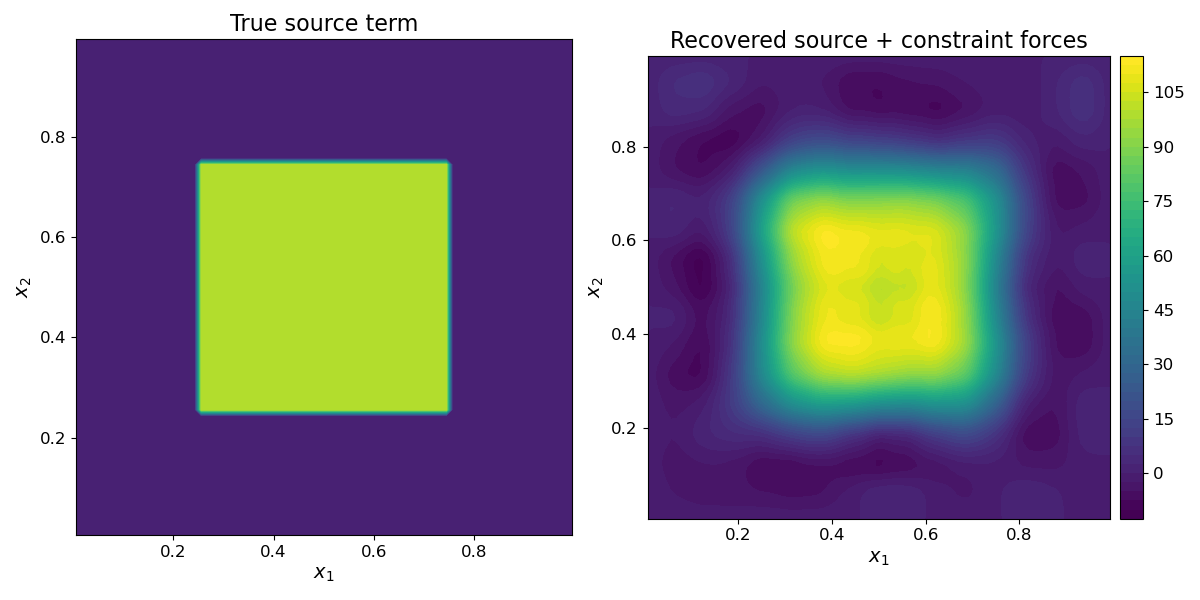}
\caption{The sum of the recovered source term and the constraint forces accurately approximates the source term that generated the data.}
\label{recovery}
\end{figure}

It is not surprising that the relative error is smaller given that the constraint forces activate to correct model misspecification. Qualitatively, the recovered source term is equivalent to that of the standard inverse problem. The lack of sharp boundaries is an inevitable consequence of the low-frequency discretization. The relative difference between the two sets of recovered parameters is given by 

\begin{equation*}
    \frac{\lVert \boldsymbol \epsilon^{\text{ECFM}} - \boldsymbol \epsilon^{\text{INV}} \rVert}{\lVert \boldsymbol \epsilon^{\text{INV}} \rVert} = 3.2 \times 10^{-2},
\end{equation*}

\noindent indicating only subtle differences in the recovered source term. However, an interesting feature of ECFM is that the constraint forces provide insight into the missing physics. Figure \ref{recovery} shows the sum of the constraint force distribution and the recovered source term. In this case, we have that 

\begin{equation*}
    s( \mathbf x) \approx b(\mathbf x; \boldsymbol \epsilon) + \sum_{i=1}^C \lambda_i \Gamma( \mathbf x - \mathbf x_i),
\end{equation*}

\noindent when $\boldsymbol \epsilon$ and $\boldsymbol \lambda$ are determined through a solution to Eq. \eqref{ecfm_noise}. The constraint force distribution acts as a discretization of the missing physics, as it can be interpreted directly as a PDE residual. As shown in \cite{rowan_physics-informed_2025}, the clear-cut physical meaning of the constraint force distribution facilitates post-processing to recover physics absent from the original parameterized model. 

\subsection{Stochastic model}

\paragraph{} We now discuss the case where the measurements are noiseless, but a component of the parameterized model is itself random. Situations of this sort arise when there are variables in an experiment which cannot be controlled, for example the case where multiple parts with different material properties are required to recover parameters of interest. To illustrate this, we consider the case of an axially-loaded Euler-Bernoulli beam with random material properties. As discussed in \cite{bauchau_structural_2009}, the governing equation for the axially-loaded beam is 

\begin{equation}
    \pdd{}{x}\qty( H(x;\omega) \pdd{u}{x}) + \beta \pdd{u}{x} = p(x), \quad x\in[0,1],
\end{equation}

\noindent where $H(x;\omega)$ is the random bending stiffness, $\beta$ is the axial load, and $p(x)$ is a distributed bending load. We take the beam to be simply supported on the left with zero rotation and shear load on the right:

\begin{equation*}
    u(0) = \pdd{u}{x}(0)=\pd{u}{x}(1) = \pd{}{x}\qty( H(1,\omega) \pdd{u}{x}(1)) = 0.
\end{equation*}

The parameterized system obeys the same governing equation and the same boundary conditions, but with the axial load $\beta$ parameterized by $\epsilon$, which is to be determined from the inverse problem. We compute the weak form of the parameterized system to obtain a numerical solution. To facilitate dealing with boundary terms from integration by parts, we do this before discretizing the trial and test space. Integrating the parameterized governing equation against an arbitrary test function $\delta w$, we obtain

\begin{equation*}
    \int_0^1 \qty(\pdd{}{x}\qty( H(x;\omega) \pdd{w}{x}) \delta w + \epsilon \pdd{w}{x} \delta w - p(x) \delta w ) dx= 0.
\end{equation*}

Note that the test function satisfies $\delta w(0)=0$ and the solution satisfies $w_{xx}(0)=0$. Integrating by parts the two terms involving spatial derivatives one time, we obtain

\begin{equation*}
    \int_0^1 \qty(\pd{}{x}\qty( H(x;\omega) \pdd{w}{x}) \pd{\delta w}{x} + \epsilon \pd{w}{x} \pd{\delta w}{x} + p(x) \delta w ) dx - \qty( \pd{}{x}\qty(H(x;\omega)\pdd{w}{x}) \delta w + \epsilon \pd{w}{x} \delta w) \Bigg |_0^1= 0.
\end{equation*}

Both boundary terms from integration by parts disappear as a consequence of the boundary conditions. Canceling the boundary terms and integrating by parts the first term a second time, the weak form of the parameterized system is 

\begin{equation*}
    \int_0^1 \qty( H(x;\omega) \pdd{w}{x} \pdd{\delta w}{x} - \epsilon \pd{w}{x} \pd{\delta w}{x} - p(x) \delta w ) dx -  H(1;\omega) \pdd{w(1)}{x} \pd{\delta w(1)}{x} = 0.
\end{equation*}

 For simplicity, we assume that the bending stiffness at the right end is not random, i.e., $H(1,\omega)=H_0$. Furthermore, we assume that the spatial basis functions satisfy all four boundary conditions by construction. We note that the remaining boundary term is a consequence of the non-zero bending moment on the right end of the beam. Plugging in the discretization $\hat w(x;\boldsymbol \theta) = \sum_{i=1}^N \theta_i f_i(x)$ and discretizing the test function with the same basis, we obtain the discretized weak form of the system. In this example, we do not perform a comparison of ECFM with the standard approach to the inverse problem. Our goal is primarily to illustrate an extension of ECFM to inverse problems with stochastic models. Thus, per the ECFM approach, we add constraint forces centered on the measurement positions. The discretized governing equation is thus

\begin{equation}\label{random_govern}
    ( \mathbf{K}^{\text{B}}(\omega) - \epsilon \mathbf K^{\text{G}} - \mathbf{K}^{\text{BC}})\boldsymbol \theta - \mathbf F - \boldsymbol \Gamma \boldsymbol \lambda = \mathbf 0,
\end{equation}

\noindent where $\mathbf K^{\text{B}}$ is the bending stiffness matrix, $\mathbf K^{\text{G}}$ is the ``geometric'' stiffness matrix corresponding to the axial load, and $\mathbf K^{\text{BC}}$ is the contribution to the stiffness matrix from the bending moment on the right end. We remark that when the axial load is sufficiently large, the effective stiffness becomes singular and the beam buckles. We ensure that the true axial load $\beta$ is such that this does not occur.

\paragraph{} Now, we assume that, unlike past numerical experiments, a series of measurements are taken at each point in space, corresponding to multiple realizations of the random variable $\omega$. We seek a formulation of the inverse problem that respects the stochastic nature of the parameterized model, which requires propagating uncertainty in the bending stiffness through the governing equation of Eq. \eqref{random_govern}. To accomplish this, we discretize the dependence of the solution parameters $\boldsymbol \theta$ on the random variable $\omega$ with a polynomial chaos expansion (PCE) \cite{ghanem_stochastic_1991}. Defining stochastic shape functions $ \{ \Psi_i(\omega) \}_{i=1}^M$, we discretize the solution parameters in the stochastic space with

\begin{equation*}
    \theta_i(\omega) = \Theta_{ij} \Psi_j(\omega), \quad j=1,2,\dots,M,
\end{equation*}

\noindent where $M$ is the number of stochastic shape functions in the discretization and $\Theta_{ij}$ is the set of unknown coefficients governing the system's response in both the physical and stochastic space. Taking the random variable governing the bending stiffness to be uniformly distributed $\omega \sim \mathcal U(0,1)$ and plugging in the discretization, the stochastic Galerkin projection of Eq. \eqref{random_govern} is obtained by multiplying by $\Psi_{\ell}(\omega)$ and computing an expectation:

\begin{equation}\label{pce}
     \Theta_{jk}\int_0^1  (K_{ij}^{\text{B}}(\omega) - \epsilon K_{ij}^{\text{G}} - K_{ij}^{\text{BC}})\Psi_k(\omega) \Psi_{\ell}(\omega)d\omega - \int_0^1 ( F_i - \Gamma_{ij} \lambda_j) \Psi_{\ell}(\omega) d\omega = 0.
\end{equation}

With the coefficients $\boldsymbol \Theta$ determined through Eq. \eqref{pce}, we use the PCE to compute statistics of the model predictions. The mean and variance of the model predictions are given by

\begin{equation}\label{stats}
    \mu_i = \mathbb E( v_i(\omega)) = \mathcal M_{ij} \Theta_{jk} \int_0^1 \Psi_k(\omega) d\omega, \quad \sigma^2_i = \mathbb E(v_i^2(\omega)) - \mu_i^2 = \mathcal M_{ij} \Theta_{jk} \mathcal{M}_{ip} \Theta_{pq} \int_0^1 \Psi_k(\omega) \Psi_q(\omega) d\omega - \mu_i^2 ,
\end{equation}

\noindent where there is no summation on the $i$ index in the expression for the variance. We decide to treat the solution parameters as an explicit function of the model parameter $\epsilon$ and the constraint forces $\boldsymbol \lambda$. In other words, $\boldsymbol \Theta = \boldsymbol \Theta(\epsilon,\boldsymbol \lambda) $ through Eq. \eqref{pce}. 

\paragraph{} At this point, we have shown how to perform a forward solve and how to use the PCE coefficients to compute the first two moments of the model prediction at the measurement points. Now, we need to compare the model predictions to the measurement data in some way, and devise a notion of optimality for the model parameter $\epsilon$ based on constraint forces. We cannot form confidence intervals on the model predictions as we did in the previous section, because their distribution is not known in advance. One strategy is to use the density of the model predictions to carry out maximum likelihood estimation (MLE) for the model parameter. Analogous to enforcing equality constraints on deterministic data, the constraint forces would then be used to enforce that the data maximizes the likelihood. This was not possible in the previous section, as only one measurement was taken at each point. In the case of a single measurement with Gaussian noise, the MLE-maximizing constraint forces would amount to an equality constraint on the data. The case of multiple measurements at each point makes likelihood maximization a possibility. However, we do not have an analytic expression for the density, as it is only defined implicitly through the measured PCE solution, given by $\mathcal M_{ij} \Theta_{jk} (\epsilon, \boldsymbol \lambda)\Psi_k(\omega)$. In fact, the density of the measurements is not even well-defined, as it is obtained through a transformation $\mathbb R^1 \rightarrow \mathbb R^C$, meaning that the pushforward measure is singular in $\mathbb{R}^C$ and no density exists in the usual sense. Workarounds exist, for example building the likelihood function as a product of the marginal densities of each component of the measurement vector \cite{lindsay_composite_1988}. Though feasible in principle, this requires pushing the probability density through the PCE expansion using the change of variables formula. When doing MLE, we then require gradients not just of the parameters of the density but of the density function itself. Obtaining the density and its gradient through the PCE expansion is intractably complex. To avoid this issue, we instead maximize a so-called ``pseudo-likelihood,'' which is a product of surrogate models for the unknown marginal density functions \cite{besag_statistical_1975}. In particular, we perform MLE with independent univariate normal distributions at each measurement point. The mean and variance of each normal distribution are given by the outputs of the PCE model through Eq. \eqref{stats}. Though an obvious approximation of the true density, this approach is principled, as a Gaussian is known to be the entropy-maximizing distribution with given first and second moments \cite{cover_elements_2006}. We define the pseudo-likelihood function as

\begin{equation*}
    \mathcal L = \prod_{i=1}^C \prod_{j=1}^D p(v_{ij} | \mu_i,\sigma_i^2),
\end{equation*}

\noindent where $D$ is the number of samples at each measurement point, $v_{ij}$ is the $j$-th replicate of the measurement at the $i$-th point, and $p(v|\mu,\sigma^2)$ is the density of a univariate normal distribution with the given mean and variance. The negative logarithm of this pseudo-likelihood function is 

\begin{equation*}
    \Lambda = - \log \mathcal L=\sum_{i=1}^C \sum_{j=1}^D \frac{1}{2} \log 2\pi \sigma_i^2 + \frac{1}{2 \sigma_i^2}( v_{ij} - \mu_i)^2.
\end{equation*}

With this pseudo-likelihood function in hand, we are now ready to introduce our approach to solving the inverse problem with the stochastic model. The constraint forces are used to maximize the pseudo-likelihood of the data, but this can be accomplished for any model parameter $\epsilon$. The solution of the inverse problem is given by the model parameter such that the likelihood-maximizing constraint forces are minimal. This optimization problem reads

\begin{equation}\label{stochastic_model}
    \begin{aligned}
        \underset{ \epsilon,\boldsymbol \lambda}{\text{argmin }} \frac{1}{2} \lVert \boldsymbol \lambda \rVert^2 \\
        \text{s.t. } \pd{\Lambda}{\boldsymbol \lambda} = \mathbf 0.
    \end{aligned}
\end{equation}

Note that the  pseudo-likelihood depends on the model parameter through the solution parameters, which determine the mean and variance of the predictions at measurement points. To clarify this point, we write

\begin{equation}\label{Lam}
    \Lambda(\epsilon, \boldsymbol \lambda) =\sum_{i=1}^C \sum_{j=1}^D  \frac{1}{2}\log 2\pi \sigma^2_i(\epsilon, \boldsymbol \lambda) + \frac{1}{2 \sigma_i^2(\epsilon, \boldsymbol \lambda )}( v_{ij} - \mu_i(\epsilon,\boldsymbol \lambda ))^2.
\end{equation}

The gradient of the negative log pseudo-likelihood with respect to the constraint forces is required for the optimization problem of Eq. \eqref{stochastic_model}:

\begin{equation}\label{gradlam}
    \pd{\Lambda}{\boldsymbol \lambda} = \sum_{i=1} ^C\sum_{j=1}^D \frac{1}{ 2\sigma_i^2} \pd{\sigma_i^2}{\boldsymbol \lambda} - (v_{ij}-\mu_i)^2 \frac{1}{2\sigma_i^4}\pd{\sigma_i^2}{\boldsymbol \lambda} - \frac{1}{\sigma_i^2}(v_{ij}-\mu_i)\pd{\mu_i}{\boldsymbol \lambda}.
\end{equation}

The gradients of the mean and variance are given through the PCE expansion as

\begin{equation}\label{sens_stats}
    \begin{aligned}
        \pd{\mu_i}{\lambda_{\ell}} = \mathcal M_{ij} \pd{\Theta_{jk}}{\lambda_{\ell}} \int_0^1 \Psi_k d\omega, \\
        \pd{\sigma_i^2}{\lambda_{\ell}} = 2 \mathcal M_{ij} \Theta_{jk} \mathcal M_{ip} \pd{\Theta_{pq}}{\lambda_{\ell}} \int_0^1 \Psi_k \Psi_q d\omega - 2 \mu_i \pd{\mu_i}{\lambda_{\ell}},
    \end{aligned}
\end{equation}

\noindent where there is again no sum over the $i$ index in the expression for the gradient of the variance. Remember that we treat the solution parameters as an explicit function of the constraint force and model parameter. Thus, the gradient of the solution parameters with respect to the constraint forces is obtained by differentiating the governing equation of Eq. \eqref{pce}. The sensitivity derivative solves the following linear system:

\begin{equation}\label{stochastic_sens}
      \pd{\Theta_{jk}}{\lambda_q}\int_0^1  (K_{ij}^{\text{B}}(\omega) - \epsilon K_{ij}^{\text{G}} - K_{ij}^{\text{BC}})\Psi_k(\omega) \Psi_{\ell}(\omega)d\omega - \int_0^1   \Gamma_{iq}  \Psi_{\ell}(\omega) d\omega = 0.
\end{equation}

The optimization problem of Eq. \eqref{stochastic_model}, as well as the gradients defined in Eqs. \eqref{gradlam}, \eqref{sens_stats}, and \eqref{stochastic_sens}, represent our proposed method for solving the inverse problem with constraint forces when the model is stochastic.

\paragraph{} We now set up a numerical experiment to test the proposed method. We take the random bending stiffness of the beam to be 

\begin{equation*}
    H(x,\omega) = H_0( 2\omega|x-1/2| - \omega), \quad \omega \sim\mathcal U(0,1).
\end{equation*}

\begin{figure}[h]
\centering
\includegraphics[width=0.99\textwidth]{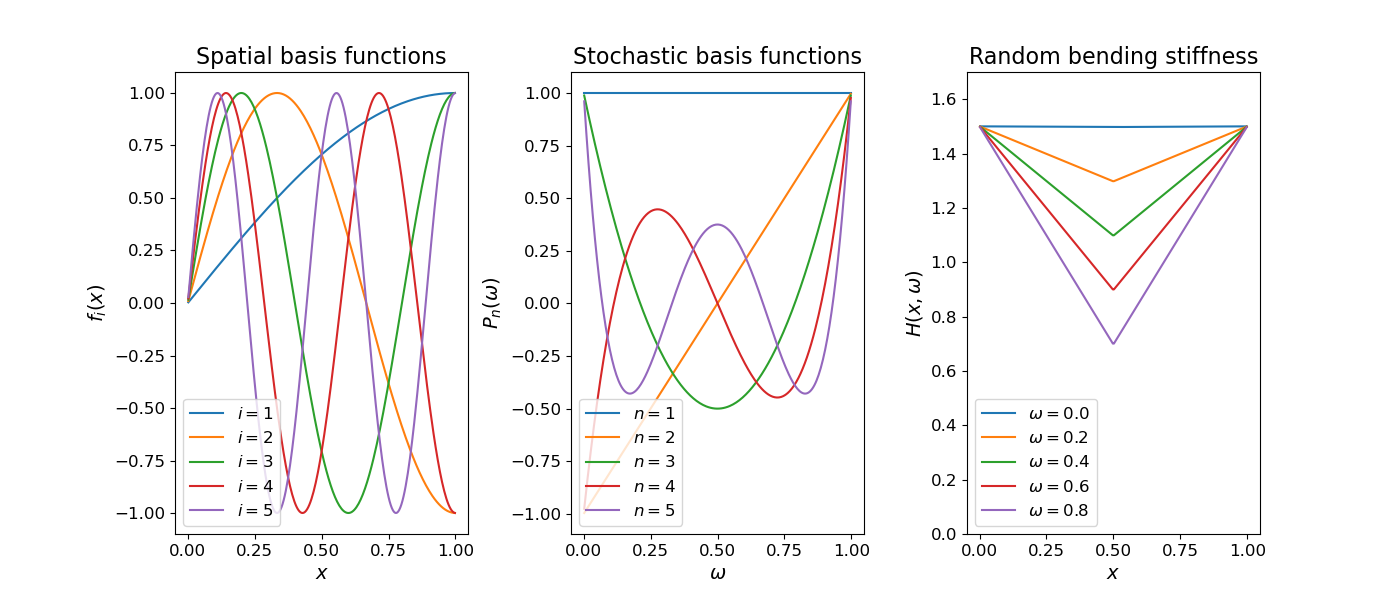}
\caption{The spatial basis functions satisfy all boundary conditions by construction (left). The basis functions for the stochastic space are orthogonal polynomials with respect to the uniform density (center). The bending stiffness distribution is an idealization of a random defect in the center of the beam, which decreases its stiffness (right).}
\label{plots}
\end{figure}

This ensures that the bending stiffness at the right endpoint is deterministic, and corresponds to a worst-case buckling load ($\omega=1$) of $\beta_{crit}=2.24$. Thus, to avoid buckling, we set the true load at $\beta=1$, and take the distributed bending load to be $p(x)=100$. The spatial shape functions are given by $f_i(x) = \sin((2i-1)\pi x/2)$, which ensure that the bending displacement satisfies all four boundary conditions by construction. The stochastic basis functions are standard Legendre polynomials shifted to the interval $\omega=[0,1]$, which we denote as $P_n(\omega)$. See Figure \ref{plots} for the two sets of basis functions and the form of the random bending stiffness. There are $C=5$ measurement points and each measurement point comprises $D=25$ replicates. Before proceeding, we verify our implementation of the PCE forward solve of Eq. \eqref{pce} and calibrate the size of the basis expansions in the physical and stochastic spaces. With zero constraint force, we solve the governing equation at discrete settings of the random variable $\omega$. In other words, we compute  

\begin{equation*}
    \boldsymbol \theta_i = ( \mathbf{K}^{\text{B}}(\omega_i) - \beta \mathbf K^{\text{G}} - \mathbf{K}^{\text{BC}})^{-1}\mathbf F 
\end{equation*}

\noindent for $\omega_i=0,1/2,1$. We then compare these predictions of the spatial form of the bending displacement to predictions obtained from the PCE solve. The reference solutions at discrete settings of the bending stiffness are computed with $15$ basis functions, and the PCE expansion uses $N=M=6$ spatial and stochastic basis functions. We remark that low-frequency solutions are expected for bending problems, given the quartic damping on high-frequency Fourier bases that arises from the fourth-order Euler-Bernoulli bending operator. We also use the PCE solution to generate data for the inverse problem. With $\epsilon=\beta=1$ and $\boldsymbol \lambda = \mathbf 0$, the measurements are governed by 

\begin{equation*}
    v_i(\omega) = \mathcal M_{ij} \Theta_{jk}(\epsilon=1, \boldsymbol \lambda = \mathbf 0)\Psi_k(\omega).
\end{equation*}

\begin{figure}[h]
\centering
\includegraphics[width=0.99\textwidth]{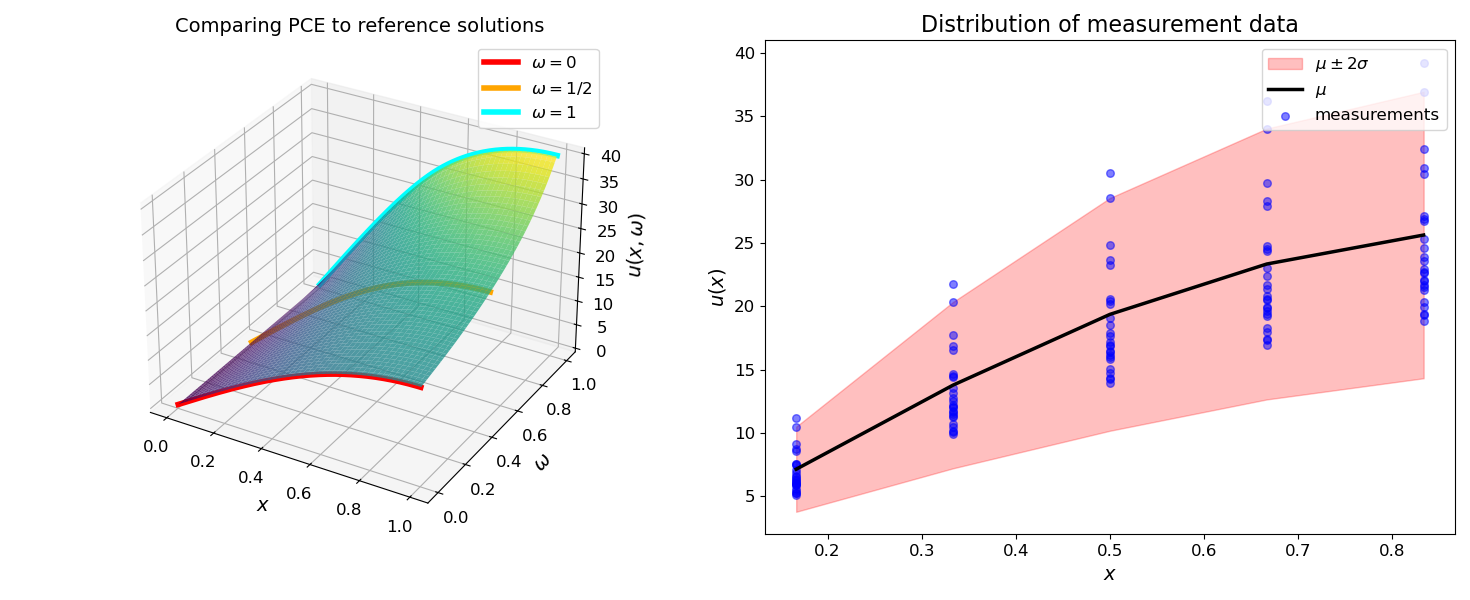}
\caption{Verifying that the PCE expansion agrees with comparatively high-fidelity solutions for the displacement at three given settings of the bending stiffness (left). The measurement data is taken at $C=5$ spatial points with $D=25$ replicates per point (right).}
\label{comparison}
\end{figure}

We generate the data $v_{ij}$ by taking $D$ random samples from $\omega \sim \mathcal U(0,1)$. The results of the comparison and the generation of measurement data are shown in Figure \ref{comparison}. The PCE solution agrees satisfactorily with the reference solutions generated by the higher-fidelity discretization. The empirical distribution of the measurement data suggests that the Gaussian surrogate used in the likelihood maximization is incorrect, as the data is skewed toward lower values of the displacement. We flag this as one potential source of error in the recovery of the model parameter.

\paragraph{} Having established the form of the discretization and generated the data, we are prepared to solve the inverse problem with the stochastic model. Given the complexity of the equality constraint on the stationarity of the pseudo-likelihood in Eq. \eqref{stochastic_model}, we choose to approximate a solution to this constrained optimization problem with a penalty approach. The problem is written in an unconstrained form as 

\begin{equation*}
    \underset{\epsilon, \boldsymbol \lambda}{\text{argmin }} \log(\frac{1}{2} \lVert \boldsymbol \lambda \rVert^2) + \alpha \Lambda( \epsilon, \boldsymbol \lambda),
\end{equation*}

\noindent where $\alpha$ is a penalty parameter and we take the logarithm of the constraint force in order to better balance the two terms in the loss. Using the scipy SLSQP optimizer with a penalty parameter of $\alpha=100$, we recover an axial load of $\epsilon=0.987$. The expected error with the true solution is given by 

\begin{equation*}
    \mathcal E = \frac{\int_0^1 \int_0^1 | u(x,\omega) - \hat w(x, \omega)| dx d \omega}{\int_0^1 \int_0^1 | u(x,\omega)| dx d \omega} = 7 \times 10^{-3}.
\end{equation*}

As expected for a problem where the parameterized model is consistent with the measurement data, the constraint force magnitude is approximately zero at the optimal setting of the model parameter. This example represents a proof of concept for extending ECFM to inverse and solution reconstruction problems in which the parameterized model contains stochastic components.


\section{The cases of parametric boundary conditions and domain geometry}

\paragraph{} In this section, we show that it is straightforward to extend ECFM to inverse or solution reconstruction problems for which the boundary conditions and/or the domain geometry are parameterized. For both these cases, we show that when the problem is discretized, the model parameters $\boldsymbol \epsilon$ appear in the discretized differential operator and/or the force vector, meaning that these cases require no method development beyond that of the original ECFM paper \cite{rowan_physics-informed_2025}. Accordingly, we derive the corresponding ECFM optimization problems but do not present numerical examples, as, at the highest level, the cases of parametric boundary conditions and geometry amount to problems that have been studied previously.

\paragraph{} To make the following exposition more concrete, we consider the case of steady-state heat conduction. First, we discuss parametric boundary conditions. Solving inverse problems for unknown or partially known boundary conditions is of interest in power electronics, where invasive measurements would be required to obtain full knowledge of thermal boundary conditions \cite{cai_physics-informed_2021}. In the case of parametric boundaries, the governing equation is

\begin{equation*}
    \begin{aligned}
        \nabla^2 w + s(\mathbf{x}) = 0 , \quad \mathbf x \in \Omega,\\
        w(\mathbf{x}) = g(\mathbf x;\boldsymbol \epsilon), \quad \mathbf x \in \partial \Omega^D, \\
        \nabla w \cdot \mathbf n = h(\mathbf x; \boldsymbol \epsilon), \quad \mathbf x \in \partial \Omega^N,
    \end{aligned}
\end{equation*}

\noindent where superscripts $D$ and $N$ denote Dirichlet and Neumann boundaries respectively. First, as is standard in the finite element and PINNs communities \cite{hughes_thomas_j_r_finite_2000, sukumar_exact_2022}, we build the inhomogeneous Dirichlet boundary condition into the discretization:

\begin{equation*}
    \hat w(\mathbf x; \boldsymbol \theta, \boldsymbol \epsilon) = \sum_{i=1}^N \theta_i f_i(\mathbf x) + g(\mathbf x; \boldsymbol \epsilon),
\end{equation*}

\noindent where the basis functions $f_i(\mathbf x)$ are zero along the Dirichlet portion of the boundaries. With slight abuse of notation, we now take $g(\mathbf x; \boldsymbol \epsilon)$ as a function which satisfies the inhomogeneous Dirichlet on the relevant portion of the boundary, but is also defined over the entire domain. Now, compute the Galerkin weak form of the governing equation and plug in the discretization:

\begin{equation*}
    \int_{\Omega}\nabla^2\qty(\sum_{i=1}^N \theta_i f_i(\mathbf x) + g(\mathbf x; \boldsymbol \epsilon)) f_j + sf_j d\Omega = 0.
\end{equation*}

Integrating by parts, we obtain

\begin{equation*}
    \int_{\Omega} \nabla \qty(\sum_{i=1}^N \theta_i f_i(\mathbf x) + g(\mathbf x; \boldsymbol \epsilon)) \cdot \nabla f_j - s f_j d\Omega - \int_{\partial \Omega^N} \nabla \qty(\sum_{i=1}^N \theta_i f_i(\mathbf x) + g(\mathbf x; \boldsymbol \epsilon)) \cdot \mathbf n f_j dS = 0.
\end{equation*}

Recognizing the boundary integral term as the Neumann boundary condition and re-arranging, we obtain a linear system for the solution parameters $\boldsymbol \theta$:

\begin{equation*}
    \qty( \int_{\Omega} \nabla f_i \cdot \nabla f_j d\Omega ) \theta_i = \int sf_j d\Omega + \int_{\partial \Omega^N} h(\mathbf x ; \boldsymbol \epsilon) f_j dS - \int_{\Omega} \nabla g(\mathbf x; \boldsymbol \epsilon) \cdot \nabla f_j d\Omega.
\end{equation*}

In symbolic notation, this expression is written as $\mathbf K \boldsymbol \theta - \mathbf F(\boldsymbol \epsilon) = \mathbf 0$. The ECFM problem to estimate the model parameters is 

\begin{equation*}
    \begin{aligned}
        \underset{\boldsymbol \epsilon, \boldsymbol \lambda, \boldsymbol \theta}{\text{argmin }} \frac{1}{2} \lVert \boldsymbol \lambda \rVert^2 \\
        \text{s.t. } \mathbf K \boldsymbol \theta - \mathbf F(\boldsymbol \epsilon) - \boldsymbol \Gamma \boldsymbol \lambda = \mathbf 0, \quad \mathcal M \boldsymbol \theta - \mathbf v = \mathbf 0.
    \end{aligned}
\end{equation*}

The case of the parameterized source term in the discretized governing equation has been treated in \cite{rowan_physics-informed_2025}. Thus, we claim that handling parameterized boundary conditions using ECFM for either solution reconstruction or an inverse problem does not require any new method development.

\paragraph{} Parametric domain geometries arise when measurements of a system's response are used to estimate the position of the boundary. A famous example of this is the question: ``can you hear the shape of a drum?'', which inquires whether the shape of a drum head uniquely determines the natural frequencies of vibration \cite{kac_can_1966}. When the domain geometry is parameterized, we write the governing equation for steady-state heat conduction as

\begin{equation*}
    \begin{aligned}
        \nabla^2 w + s(\mathbf{x}) = 0 , \quad \mathbf x \in \Omega(\boldsymbol \epsilon),\\
        w(\mathbf{x}) = 0, \quad \mathbf x \in \partial  \Omega(\boldsymbol \epsilon),
    \end{aligned}
\end{equation*}

The assumption of homogeneous Dirichlet boundaries is made for simplicity, but places no meaningful restrictions on our conclusions. We discretize the solution with $\hat w(\mathbf x; \boldsymbol \theta;\boldsymbol \epsilon) = \sum_{i=1}^N \theta_i f_i(\mathbf x ; \boldsymbol \epsilon)$, where the basis functions are taken to depend on the geometry parameters by virtue of satisfying the homogeneous Dirichlet boundary conditions by construction. For example, with a finite element basis, this dependence on the model parameters is through the nodes in the mesh defining the computational domain. Computing the Galerkin weak form of the governing equation and plugging in the discretization, we obtain

\begin{equation*}
    \qty(\int_{\Omega(\boldsymbol \epsilon)} \nabla f_i \cdot \nabla f_j  d\Omega) \theta_i  -\int_{\Omega(\boldsymbol \epsilon)} s f_j d\Omega = 0.
\end{equation*}

Both the basis functions and the integration domain depend on the model parameters. However, this equation can be written simply as $\mathbf K(\boldsymbol \epsilon) \boldsymbol \theta - \mathbf F(\boldsymbol \epsilon) = \mathbf 0$. The corresponding ECFM problem is

\begin{equation*}
    \begin{aligned}
        \underset{\boldsymbol \epsilon, \boldsymbol \lambda, \boldsymbol \theta}{\text{argmin }} \frac{1}{2} \lVert \boldsymbol \lambda \rVert^2 \\
        \text{s.t. } \mathbf K(\boldsymbol \epsilon) \boldsymbol \theta - \mathbf F(\boldsymbol \epsilon) - \boldsymbol \Gamma \boldsymbol \lambda = \mathbf 0, \quad \mathcal M \boldsymbol \theta - \mathbf v = \mathbf 0.
    \end{aligned}
\end{equation*}

Here, both the force vector and the stiffness matrix depend on the model parameters. This is the same kind of problem that is obtained when both the source term and differential operator are parameterized. Parameterized sources and differential operators have been treated separately in \cite{rowan_physics-informed_2025}, but we argue that no new developments are required to treat the case where both are parameterized simultaneously. This establishes that ECFM can treat parameterized domain geometries without any fundamental modifications.

\paragraph{} Though the above discussion was carried out in the context of steady-state heat conduction, we claim that the conclusions generalize to nonlinear and/or vector-valued PDEs. Whether in the setting of an inverse problem or solution reconstruction, parameters of the boundary conditions appear exclusively in the force vector when the problem is discretized. Similarly, when the domain geometry is to be recovered from measurements, the model parameters show up in both the discretized differential operator and the force vector. Although this type of parametric dependence is more familiar from volumetric source terms or material coefficients, from the perspective of the discretized governing equation, estimating boundary conditions or domain geometry is an equivalent problem.


\section{Conclusion}

\paragraph{} The purpose of this work was two-fold: to demonstrate the use of ECFM for inverse problems, and to discuss extensions of the method to novel problem types. In particular, we discussed extensions to dynamic problems, stochastic models, parametric boundary conditions and domain geometries, and also introduced a streamlined approach to handle noisy measurements. The cases of parametric boundary conditions and domain geometries were the most straightforward, as they were shown to follow from existing test cases of ECFM. In treating the dynamic problem, we introduced time-dependent constraint forces at measurement points, and performed sensitivity analyses to remove constraints from the optimization problem. In this noiseless setting, both the standard formulation of the inverse problem and ECFM recovered the exact solution within the allotted $250$ optimization steps. Then, noting the prohibitively intrusive implementation of the approach for handling noisy measurement data introduced in a previous work, we devised a new method for measurements with a known noise distribution. This formulation relied on choosing the model parameter such that the magnitude of the constraint force that enforced agreement between statistics of the discrepancy and the measurement noise was minimal. Though this disrupted the inner-outer loop structure that has characterized past treatments of the ECFM problem, we wrote down an inequality-constrained optimization problem which could be solved with standard optimizers. Lastly, we devised a method based on a polynomial chaos expansion and maximum likelihood estimation to solve inverse problems with constraint forces when the parameterized governing equation was stochastic. 

\paragraph{} We reiterate a point made previously---whereas a standard inverse problem minimizes the model's discrepancy with the data in the sense of ``displacement,'' the constraint force formulation does so in the sense of ``force.'' This terminology takes inspiration from solid mechanics. More generally, a standard inverse problem operates on the system state whereas the ECFM works with the source terms. We believe the fact that the ECFM objective more heavily weights errors in stiffer parts of the domain to be philosophically interesting. Is it clear that the standard inverse problem's equal treatment of all errors is justified? Should errors which require larger forces to correct not be prioritized? We have no interest in arguing that one or another method is ``better.'' This work is simply an exploration of an alternative formulation of deterministic inverse problems. In the case that the parameterized model is consistent with the data, the zero constraint force and zero error solution agree on the recovered model parameters. In the case where the parameterized model is inconsistent---which is really the more interesting case---it is not clear how to compare the performance of the two methods. In general, we do not expect to recover the true model parameters with either method, nor the exact solution field. As our second numerical example shows, one unique feature of ECFM is that it provides a discretization of the missing physics.

\paragraph{} At this point, the constraint force method has been studied in the context of solution reconstruction and/or inverse problems for 1) parameterized differential operators, 2) parameterized source terms, 3) noisy measurements (two strategies), 4) dynamic problems, and 5) stochastic models. We have shown that parametric boundary conditions and domain geometries are special cases of the first two problem types once the PDE is discretized. Many more example problems will need to be tried to establish the efficacy of this method compared to the standard approach to deterministic inverse problems. Our purposes here are to lay the groundwork for such studies. Future work will focus on exploring the explicit constraint force method in another novel context: optimal experimental design.



\end{document}